\documentclass[12pt]{article}
\usepackage{theorem}
\usepackage{amssymb,amsmath,mathrsfs}
\newtheorem{prop}{}[section]

{\theorembodyfont{\upshape} \newtheorem{rema}[prop]{}}
\newcommand{\boma}[1]{{\mbox{\boldmath $#1$} }}
\begin{document}
\newcommand{\mper}[1]{\stackrel{\mbox{\footnotesize{(#1)}}}{\leqs} }
\def\gammaa{\beta}
\def\epp{\delta}
\def\bet{\delta}
\def\gaa{\beta}
\def\mpe{\barray{ccc} ~ \\ \leqs \\ \stackrel{(\Co_{2 3 d})}{\stackrel{~}{\stackrel{~}{~}}} \farray}
\def\Co{{\mathfrak S}}
\def\Bst{B \!\raisebox{-0.1cm}{$\scriptstyle{s t}$}}
\def\hs{\hspace{-0.5cm}}
\def\ugp{~\overline{\thicksim}~}
\def\ugm{~\underline{\thicksim}~}
\def\n{t}
\def\S{S}
\def\ppsi{\vartheta}
\def\PPsi{\Theta}
\def\L{M}
\def\M{N}
\def\N{V}
\def\R{R}
\def\half{{1 \over 2}}
\def\II{M}
\def\JJ{N}
\def\ga{\gamma}
\def\mmu{\nu}
\def\nnu{\mu}
\def\F{G}
\def\f{g}
\def\Fs{\mathscr{G}}
\def\Kap{\mathscr{K}}
\def\STB{\scriptsize B}
\def\STBB{\scriptsize BB}
\def\STF{\scriptsize F}
\def\STFF{\scriptsize FF}
\def\TB{\scriptsize(B)}
\def\TBB{\scriptsize(BB)}
\def\TF{\scriptsize(F)}
\def\TFF{\scriptsize(FF)}
\def\p{\alpha}
\def\pp{\beta}
\def\XXX{\mathscr X}
\def\Piu{\mathscr P}
\def\Men{\mathscr N}
\def\ffi{\varphi}
\def\ES{{\mathcal S}}
\def\J{{\mathscr J}}
\def\H{{\mathscr H}}
\def\K{{\mathscr K}}
\def\Kp{{\mathscr K}'}
\def\kp{k'}
\def\scrscr{\scriptscriptstyle}
\def\scr{\scriptstyle}
\def\dd{\displaystyle}
\def\z{w}
\def\x{w}
\def\w{\eta}
\def\B{ B_{\mbox{\scriptsize{\textbf{C}}}} }
\def\Bc{ \overline{B}_{\mbox{\scriptsize{\textbf{C}}}} }
\def\ppartial{\overline{\partial}}
\def\d{\hat{d}}
\def\TT{T}
\def\Hinf{ H^{\infty}(\reali^d, \complessi) }
\def\Hn{ H^{n}(\reali^d, \complessi) }
\def\Hm{ H^{m}(\reali^d, \complessi) }
\def\Hell{ H^{\ell}(\reali^d, \complessi) }
\def\Ht{ H^{t}(\reali^d, \complessi) }
\def\Ha{ H^{\d}(\reali^d, \complessi) }
\def\Ld{L^{2}(\reali^d, \complessi)}
\def\Lpi{L^{p}(\reali^d, \complessi)}
\def\Lq{L^{q}(\reali^d, \complessi)}
\def\Lr{L^{r}(\reali^d, \complessi)}
\def\Knb{K^{best}_n}
\def\k{\mbox{{\tt k}}}
\def\D{\mbox{{\tt D}}}
\def\g{ {\textbf g} }
\def\QQQ{ {\textbf Q} }
\def\AAA{ {\textbf A} }
\def\gr{\mbox{graph}~}
\def\Q{$\mbox{Q}_a$~}
\def\PZ{$\mbox{P}^{0}_a$~}
\def\PZAL{$\mbox{P}^{0}_\alpha$~}
\def\PL{$\mbox{P}^{1/2}_a$~}
\def\PU{$\mbox{P}^{1}_a$~}
\def\PK{$\mbox{P}^{k}_a$~}
\def\PKU{$\mbox{P}^{k+1}_a$~}
\def\PI{$\mbox{P}^{i}_a$~}
\def\Pell{$\mbox{P}^{\ell}_a$~}
\def\PTM{$\mbox{P}^{3/2}_a$~}
\def\AZ{$\mbox{A}^{0}_r$~}
\def\AU{$\mbox{A}^{1}$~}
\def\epsilona{\epsilon^{\scriptscriptstyle{<}}}
\def\epsilonb{\epsilon^{\scriptscriptstyle{>}}}
\def\lgraffa{ \mbox{\Large $\{$ } \hskip -0.2cm}
\def\rgraffa{ \mbox{\Large $\}$ } }
\def\restriction{ \stackrel{\setminus}{~}\!\!\!\!|~}
\def\m{m}
\def\Fre{Fr\'echet~}
\def\I{{\mathcal N}}
\def\ap{{\scriptscriptstyle{ap}}}
\def\fiap{\varphi_{\ap}}
\def\BBB{ {\textbf B} }
\def\EEE{ {\textbf E} }
\def\FFF{ {\textbf F} }
\def\TTT{ {\textbf T} }
\def\KKK{ {\textbf K} }
\def\FFi{ {\bf \Phi} }
\def\GGam{ {\bf \Gamma} }
\def\a{a}
\def\ep{\epsilon}
\def\parn{\par\noindent}
\def\teta{M}
\def\elle{L}
\def\ro{\rho}
\def\al{\alpha}
\def\si{\sigma}
\def\be{\beta}
\def\de{\delta}
\def\la{\lambda}
\def\te{\vartheta}
\def\ch{\chi}
\def\complessi{{\textbf C}}
\def\reali{{\textbf R}}
\def\interi{{\textbf Z}}
\def\naturali{{\textbf N}}
\def\bT{{\textbf T}}
\def\T1{{\textbf T}^{1}}
\def\EE{{\mathcal E}}
\def\FF{{\mathcal F}}
\def\EFFE{{\mathscr F}}
\def\GG{{\mathcal C}}
\def\PP{{\mathcal P}}
\def\QQ{{\mathcal Q}}
\def\Np{{\hat{N}}}
\def\Lp{{\hat{L}}}
\def\Jp{{\hat{J}}}
\def\Pp{{\hat{P}}}
\def\Pip{{\hat{\Pi}}}
\def\Vp{{\hat{V}}}
\def\Ep{{\hat{E}}}
\def\Fp{{\hat{F}}}
\def\Gp{{\hat{G}}}
\def\Ip{{\hat{I}}}
\def\Tp{{\hat{T}}}
\def\Mp{{\hat{M}}}
\def\La{\Lambda}
\def\Ga{\Gamma}
\def\Si{\Sigma}
\def\Upsi{\Upsilon}
\def\Gag{{\check{\Gamma}}}
\def\Lap{{\hat{\Lambda}}}
\def\Sip{{\hat{\Sigma}}}
\def\Upsig{{\check{\Upsilon}}}
\def\Kg{{\check{K}}}
\def\ellp{{\hat{\ell}}}
\def\j{j}
\def\jp{{\hat{j}}}
\def\Stir{{\mathscr S}}
\def\Ii{{\mathscr H}}
\def\BB{{\mathscr B}}
\def\LL{{\mathcal L}}
\def\SS{{\mathscr S}}
\def\DD{{\mathcal D}}
\def\VV{{\mathcal V}}
\def\WW{{\mathcal W}}
\def\OO{{\mathcal O}}
\def\RR{{\mathcal R}}
\def\AA{{\mathscr A}}
\def\CC{{\mathscr C}}
\def\NN{{\mathcal N}}
\def\WW{{\mathcal W}}
\def\HH{{\mathcal H}}
\def\XX{{\mathcal X}}
\def\YY{{\mathcal Y}}
\def\ZZ{{\mathcal Z}}
\def\UU{{\mathcal U}}
\def\XX{{\mathcal X}}
\def\RR{{\mathcal R}}
\def\cir{{\scriptscriptstyle \circ}}
\def\circa{\thickapprox}
\def\vain{\rightarrow}
\def\leqs{\leqslant}
\def\geqs{\geqslant}
\def\ss{s}
\def\vains{\stackrel{\ss}{\rightarrow}}
\def\parn{\par \noindent}
\def\salto{\vskip 0.2truecm \noindent}
\def\spazio{\vskip 0.5truecm \noindent}
\def\vs1{\vskip 1cm \noindent}
\def\fine{\hfill $\diamond$ \vskip 0.2cm \noindent}
\newcommand{\rref}[1]{(\ref{#1})}
\def\beq{\begin{equation}}
\def\feq{\end{equation}}
\def\beqq{\begin{eqnarray}}
\def\feqq{\end{eqnarray}}
\def\barray{\begin{array}}
\def\farray{\end{array}}
%%%%%%%%% THIS NUMBERS EQUATIONS BY SECTIONS %%%%%%%%%%%%%
\makeatletter
\@addtoreset{equation}{section}
\renewcommand{\theequation}{\thesection.\arabic{equation}}
%\thesection instead of \arabic{section} for correct equation numbering
% in appendices
\makeatother
%%%%%%%%%%%%%%%%%%%%%%%%%%%INTESTAZIONE%%%%%%%%%%%%%%%%%%%%%%%%%%%%%%%
\begin{titlepage}
\begin{center}
{\huge New results on multiplication \\ in Sobolev spaces.}
\end{center}
\vspace{1truecm}
\begin{center}
{\large
Carlo Morosi${}^1$, Livio Pizzocchero${}^2$} \\
\vspace{0.5truecm}
${}^1$ Dipartimento di Matematica, Politecnico di
Milano, \\ P.za L. da Vinci 32, I-20133 Milano, Italy \\
e--mail: carlo.morosi@polimi.it \\
${}^2$ Dipartimento di Matematica, Universit\`a di Milano\\
Via C. Saldini 50, I-20133 Milano, Italy\\
and Istituto Nazionale di Fisica Nucleare, Sezione di Milano, Italy \\
e--mail: livio.pizzocchero@unimi.it
\end{center}
\vspace{1truecm}
\begin{abstract}
We consider the Sobolev (Bessel potential) spaces $H^\ell(\reali^d, \complessi)$,
and their standard norms $\|~\|_\ell$ (with $\ell$ integer or noninteger).
We are interested in the unknown sharp constant $K_{\ell m n d}$
in the inequality $\| f g \|_{\ell} \leqs K_{\ell m n d} \| f \|_{m}
\| g \|_n$ ($f \in H^m(\reali^d, \complessi)$, $g \in H^n(\reali^d, \complessi)$;
$0 \leqs \ell \leqs m \leqs n$, $m + n - \ell > d/2$); we derive upper and
lower bounds $K^{\pm}_{\ell m n d}$ for this constant.
As examples, we give a table of
these bounds for $d=1$, $d=3$ and many values of $(\ell, m, n)$; here the
ratio $K^{-}_{\ell m n d}/K^{+}_{\ell m n d}$ ranges between $0.75$ and $1$
(being often near $0.90$, or larger), a fact indicating that the bounds
are close to the sharp constant.
Finally, we discuss the asymptotic behavior of the upper and lower
bounds for $K_{\ell, b \ell, c \ell, d}$ when $1 \leqs b \leqs c$ and
$\ell \vain + \infty$. As an example, from this analysis we obtain the
$\ell \vain +\infty$ limiting behavior of the sharp constant $K_{\ell, 2 \ell, 2 \ell, d}$;
a second example concerns the $\ell \vain + \infty$ limit for $K_{\ell, 2 \ell, 3 \ell, d}$.
The present work generalizes our previous paper \cite{uno},
entirely devoted to the constant $K_{\ell m n d}$ in the special case
$\ell = m = n$; many results given therein can be recovered here
for this special case.
\end{abstract}
\vspace{1truecm} \noindent \textbf{Keywords:} Sobolev spaces,
inequalities, pointwise multiplication. \par \vspace{0.4truecm} \noindent \textbf{AMS 2000
Subject classifications:} 46E35, 26D10, 47A60.
\par
\end{titlepage}
\section{Introduction and preliminaries.} \label{intro}
The present work generalizes some results of ours \cite{uno} on pointwise multiplication in the Sobolev
(or Bessel potential) spaces $H^{\ell}(\reali^d, \complessi)$ (see the forthcoming
Eqs. \rref{incid} \rref{repfur} for a precise definition of these spaces and of their norms).
In the cited work, we derived upper and lower bounds for the sharp constant $K_{\ell d}$ in
the inequality
\beq \| f g \|_{\ell} \leqs K_{\ell d} \| f \|_{\ell} \| g \|_{\ell}
\qquad \mbox{for}~f, g \in H^{\ell}(\reali^d, \complessi), ~~\ell > d/2~. \label{inequno} \feq
Here, we derive bounds for the sharp constant $K_{\ell m n d}$ in the inequality
\beq \| f g \|_{\ell} \leqs K_{\ell m n d} \| f \|_{m} \| g \|_{n}
\qquad \mbox{for $f \in H^{m}(\reali^d, \complessi)$, $g \in H^{n}(\reali^d, \complessi)$,}
\label{inequi} \feq
$$  \ell, m, n \in \reali,~~0 \leqs \ell \leqs m \leqs n,~~n + m - \ell > d/2~; $$
this becomes \rref{inequno} for $\ell = m = n$. The relation
$H^{m}(\reali^d, \complessi) H^{n}(\reali^d, \complessi) \subset H^\ell(\reali^d, \complessi)$
and the inequality \rref{inequi} are well known for the indicated values of
$\ell,m,n$ (see e.g. \cite{Ben}, Part 5); however, to
the best of our knowledge, no quantitative
analysis seems to have been done for the related constants. \parn
One of the motivations to analyze the constants in this inequality and similar ones is the
same indicated in \cite{uno}: this analysis allows to infer \textsl{a posteriori
estimates} on the error of most approximation methods for semilinear
evolutionary PDEs with polynomial nonlinearities, and also to get
bounds on the time of existence for their exact solutions (see in particular \cite{approxi}, where
we considered a nonlinear heat equation and the Navier-Stokes equations). This is just one of the
possible applications: in fact, inequalities of the type  \rref{inequno} \rref{inequi}
and similar ones are relevant for several reasons in many areas of mathematical physics, including
the $\varphi^4$ quantum field theory and the analysis of the Lieb functional in electronic
density theory \cite{Lieb} \cite{Lam}. \parn
Let us fix the attention to \rref{inequi}. Finding the
sharp constant $K_{\ell m n d}$ is clearly difficult; for this reason, and even in view of applications to PDEs,
one can be satisfied to derive two-sided bounds
\beq K^{-}_{\ell m n d} \leqs K_{\ell m n d} \leqs K^{+}_{\ell m n d}~, \label{kmp} \feq
where the lower bound $K^{-}_{\ell m n d}$ is sufficiently close to the
upper bound $K^{+}_{\ell m n d}$: this is the same attitude proposed in \cite{uno}
for the constant $K_{\ell d}$ of \rref{inequno}.  \parn
In the present paper, we produce the following upper and lower bounds. \parn
(i) First of all, we establish what we call the ``$\SS$-function'' upper bound $K^{\SS}_{\ell m n d}$; this is
obtained maximizing a suitable function $\SS_{\ell m n d} : [0,+\infty) \vain (0,+\infty)$
(which is, up to a factor, a generalized hypergeometric function).
In the special case $\ell=0$, we derive as well a
``H\"older'' upper bound $K^{\Ii}_{0 m n d}$; this is
obtained from the H\"older and from the Sobolev imbedding inequalities. \parn
(ii) Next, we present a number of lower bounds; all of them are obtained directly from
Eq. \rref{inequi}, choosing for $f$, $g$
some convenient trial functions (generally depending on
certain parameters, to be fixed optimally). Different choices of the trial functions
yield the  so-called ``Bessel'' lower bound $K^{B_{s t}}_{\ell m n d}$,
the ``Fourier'' lower bound $K^{F}_{\ell m n d}$ and the ``S-constant'' lower bound
$K^{S}_{\ell \ell n d}$ (holding for $m = \ell$ only). \parn
The above terminology for the upper and lower bounds is used only
for convenience: the terms ''$\SS$-function'', etc., recall some distinguished
function or feature appearing in the construction of these bounds. For all $\ell, m, n, d$,
from the available upper and lower bounds one can extract the best ones, indicated
with $K^{\pm}_{\ell m n d}$: so, $K^{+}_{\ell m n d}$ is the minimum of the
upper bounds in (i)  and $K^{-}_{\ell m n d}$ is the maximum of the
lower bounds in (ii). \parn
To exemplify the above framework, the paper presents a table of
upper and lower bounds $K^{\pm}_{\ell m n d}$ in dimension $d=1$ and $d=3$,
for a set of values of $\ell, m, n$; in each case, informations are provided on the type
of bound employed, and on its practical computation. In all cases
presented in the table, the ratio $K^{-}_{\ell m n d}/K^{+}_{\ell m n d}$
ranges between $0.75$ and $1$, often reaching a value larger than $0.90$; so, our bounds
are not far from the sharp constant $K_{\ell m n d}$.
It would not be difficult to build similar tables,
for different values of $\ell, m, n$ (even non integer) and $d$. \parn
The final step in our analysis is the asymptotics of some available upper and lower bounds,
when $\ell, m, n$ go to infinity (and $d$ is fixed). This generalizes an analysis
performed in \cite{uno}, where we proved for the constant $K_{\ell d}$ in \rref{inequno} the relations
$$ 0.793 \, T_d ~{ (2/\sqrt{3})^\ell \over \ell^{d/4} } \Big[ 1 + O({1 \over \ell}) \Big] \leqs
K_{\ell d} \leqs T_d
~{ (2/\sqrt{3})^\ell \over \ell^{d/4} } \Big[ 1 + O({1 \over \ell}) \Big]
\qquad \mbox{for $\ell \vain + \infty$}~, $$
\beq T_d := {3^{d/4 + 1/4} \over 2^d \pi^{d/4} }~ \label{eqintro} \feq
(to be intended as follows: $K_{\ell d}$ has upper and lower bounds behaving like the
right and left hand side of the above equation). \parn
In the present paper, some of our bounds
on the sharp constant $K_{\ell, b \ell, c \ell, d}$ are investigated for $\ell \vain + \infty$ and fixed $b, c, d$
($1 \leqs b \leqs c$). To exemplify our results, let us report the conclusions
arising for $b=c=2$ and $b=2, c=3$, respectively. In the first case we grasp the
limiting behavior of the sharp constant, which is the following:
\beq K_{\ell, 2 \ell, 2 \ell, d} = {1 + O(1/\ell) \over (16 \pi \ell)^{d/4}}
\qquad \mbox{for $\ell \vain + \infty$}~; \label{find} \feq
the above result is inferred from the analysis of suitable upper and lower bounds for $K_{\ell,
2 \ell, 2 \ell, d}$, both of them behaving like the right hand side of \rref{find}
when $\ell \vain + \infty$. \parn
In the second case, we find
\beq  {1 + O(1/\ell) \over (23 \pi  \ell)^{d/4}} \leqs
K_{\ell, 2 \ell, 3 \ell, d} \mpe {1 + O(1/\ell) \over
(20 \pi \ell)^{d/4}} \qquad \mbox{for $\ell \vain + \infty$}~. \label{findsec} \feq
\vskip -0.4cm \noindent
The subscript $\scriptstyle{(\Co_{2 3 d})}$ in Eq. \rref{findsec} means that the indicated upper bound holds
under a certain condition $\Co_{2 3 d}$, dealing with the maximum of a
hypergeometric-like function; we have numerical indications
that the condition is satisfied for all $d$, as explained later in the paper.
\vskip 0.2cm \noindent
\textbf{Organization of the paper.} In the sequel of the present
section we fix a few notations, and review
some standard properties of the special functions
employed throughout the paper (Bessel, hypergeometric, etc.);
an integral identity about
Bessel functions presented here, and seemingly less trivial, is proved for completeness
in Appendix \ref{appef}.
Again in this section, we review the definition
of the spaces $H^\ell(\reali^d, \complessi)$. (Some facts reported in this section
were already mentioned in \cite{uno}; they have been reproduced
to avoid continuous, annoying citation of small details
from the previous work).
\parn
In Section \ref{desc} we present our upper and lower bounds
on $K_{\ell m n d}$, of all the types mentioned before (e.g., the ``$\SS$-function'' upper bound,
the ``Bessel'' lower bound, and so on); most proofs about
these bounds are given later, in Sections \ref{upper}, \ref{pbes}, \ref{secdir}.
\parn
In Section \ref{sectable} we describe the practical computation of the bounds
in Section \ref{desc}, and present the already mentioned table of upper and
lower bounds $K^{\pm}_{\ell m n d}$, for $d=1,3$ and many values of $\ell, m,n$;
further details on the construction of the table are given in Appendix
\ref{apptab}.
\parn
In Section \ref{secasimp} we describe the asymptotics of some upper and lower bounds
for $K_{\ell, b \ell, c \ell, d}$, when $1 \leqs b \leqs c$ and $\ell \vain + \infty$; as
examples we consider the cases $(b,c) = (2,2)$ and $(2,3)$, yielding the previous
mentioned results \rref{find} \rref{findsec}. Most statements of Section \ref{secasimp}
are proved in Section \ref{secbcd}. \parn
\vskip 0.2cm\noindent
\textbf{Some basic notations and facts.} Throughout the paper: \parn
(i) $\naturali$ stands for $\{0,1,2,...\}$, $\naturali_0$ means $\naturali
\setminus \{0\}$. We often consider the sets $-\naturali = \{0,-1,-2,....\}$, $2 \naturali = \{0,2,4,...\}$,
$2 \naturali + 1 = \{1,3,5,...\}$ and $\naturali + {1 \over 2} = \{ {1 \over 2}, {3 \over 2},
{5 \over 2},...\}$. \parn
(ii) We use the double factorial
\beq (-1)!! := 1~; \qquad s!! := 1 \cdot 3 \cdot .... \cdot (s - 2) s
\quad \mbox{for $s \in 2 \naturali + 1$}~. \label{semif} \feq
(iii) The Pochhammer symbol
of $a \in \reali$, $i \in \naturali$ is
\beq (a)_i := 1~\mbox{if $i = 0$}, \qquad (a)_i := a ( a + 1) ... (a + i - 1)~\mbox{if $i > 0$}; \label{poch} \feq
note that
\beq (-s)_i = 0 \qquad \mbox{for $s \in \naturali$, $i > s$}~. \feq
(iv) We work in any space dimension $d \in \naturali_0$. The standard
inner product and Euclidean norm of $\reali^d$ are denoted by
$\bullet$ and $|~|$, respectively. The running variable over $\reali^d$ is written $x = (x^1,...,x^d)$
(or $k$, when $\reali^d$ is viewed as the space of ``wave vectors'' for the Fourier transform);
the Lebesgue measure of $\reali^d$ is indicated with $d x$ (or $d k$). \parn
For future citation, we record here the familiar formula for integrals over $\reali^d$ of radially symmetric
functions; this is the equation
\beq \int_{\reali^d} \!\!\! d x~ \varphi(| x |) = {2 \, \pi^{d/2} \over \Gamma(d/2)}
\int_{0}^{+\infty} \!\!\! d r~ r^{d-1} \varphi(r)~, \label{rag} \feq
holding for all sufficiently regular real (or complex) functions $\varphi$ on $(0,+\infty)$
(when dealing with integrals on the
"wave vector" space $(\reali^d, d k)$, the radius $r$ is renamed $\rho$).
\vskip 0.2cm \noindent
\textbf{Some special functions.} The independent variables and the parameters appearing in the special
functions that we consider
are \textsl{real}, unless the use of complex numbers is explicitly declared; consequently,
the notion of analyticity often employed in relation with such functions is intended
in the real sense. We take \cite{Car} as a general reference on real analyticity; in particular, we
frequently refer to the principle of analytic continuation as stated
in Corollary 2, page 122 of the cited book. \parn
We take \cite{Abr} \cite{Luke} \cite{Olv} \cite{Wat}  as standard references for special functions.
In this paper, we frequently use: the Gamma function $\Gamma$;
the Bessel functions of the first kind $J_{\nu}$,
the modified Bessel functions of the first kind
$I_{\nu}$ and the modified Bessel functions of the second
kind, or Macdonald functions, $K_\nu$; the generalized hypergeometric functions
${}_p F_q$, especially in the cases $p=2,q=1$ (the usual Gaussian hypergeometric
function) and $p=3, q = 2$. \parn
Concerning the Gamma function, we often use:
the integral representation
\beq \Gamma(\al) = \int_{0}^{+\infty} \!\!\!\! \!\! \!\! d p \, p^{\al - 1} e^{-p}~\qquad \mbox{for $\al \in (0,+\infty)$}~,
\label{igamma} \feq
the elementary relations
\beq \Gamma(k + 1)= k!~, \qquad
\Gamma(\alpha + k) = (\alpha)_k \Gamma(\alpha) \qquad \mbox{for $k \in \naturali$}~, \label{elem} \feq
the duplication formula
\beq \Gamma(2 \al) = {2^{2 \al - 1} \over \sqrt{\pi}}~\Gamma(\al+{1 \over 2}) \Gamma(\al) \label{dupl}~, \feq
the integral identity
\beq \int_{0}^{1} \!\!\!
d t~{t^{\alpha-1} (1 - t)^{\beta-1}} = {\Gamma(\al) \Gamma(\be) \over \Gamma(\al + \be)}
\qquad \mbox{for $\al, \be \in (0,+\infty)$}, \label{recall} \feq
and the asymptotics
\beq {\Gamma(\al + \mu) \over \Gamma(\al + \nu)} = \al^{\mu-\nu} [1 + O({1 \over \alpha})] \qquad \mbox{for
$\mu,\nu \in \reali$, $\al \vain + \infty$}~. \label{rapgamma} \feq
As for the Macdonald functions, we recall that
\beq K_{\nu}(w) =
\sqrt{\pi}~{e^{-w}}~
\sum_{i=0}^{\nu -1/2} {( 2 \nu - i - 1)! \over i! (\nu - i - 1/2)!}~ {1 \over (2 w)^{\nu-i}}
\quad \mbox{for}~\nu \in \naturali + {1 \over 2},~w \in \reali~.
\label{cita} \feq
The list of results we need about ${}_{p} F_{q}$ functions is longer, and
wholly occupies the next paragraph. \parn
\textbf{On (generalized) hypergeometric functions.} Most of the facts reported hereafter
on the ${}_{p} F_{q}$ hypergeometric functions are derived from \cite{Luke};
we will occasionally mention other references. Let
\beq p, q \in \naturali~, \qquad \alpha_1, ..., \alpha_p \in \reali, ~~~\bet_1, ..., \bet_q
\in \reali\setminus (-\naturali)~; \feq
for $k=0,1,2,...$
we associate to the parameters $\alpha_1,...,\bet_q$ the Pochhammer's symbols
$(\alpha_1)_k$, $... ,(\alpha_p)_k$, $(\bet_1)_k,$ $... ,$ $(\bet_q)_k$, noting
that $(\bet_i)_k \neq 0$ due to the assumptions on $\bet_i$. If
$\z$ is a real variable, the standard definition
\beq {~}_{p} F_{q}(\alpha_1, ... ,\alpha_p; \bet_1,...,\bet_q; \z)
:= \sum_{k=0}^{+\infty} {(\alpha_1)_{k} ... (\alpha_p)_k \over (\bet_1)_k...(\bet_q)_k}
\, {\z^{k} \over k!}  \label{partpq} \feq
makes sense when the above power series in $\z$ converges; this happens, in particular, if
\beq p=q~, \qquad \z \in \reali \feq
or
\beq p = q+1~, \qquad \z \in (-1,1)~, \label{limiz} \feq
or
\beq p,q~ \mbox{arbitrary},~\alpha_i = - \ell~ \mbox{for some $i \in \{1,...,p\}$ and $\ell
\in \naturali$}, ~\z \in \reali~; \label{limizz} \feq
in the third case we have
$(\alpha_i)_k =0$ for $k > \ell$, so the series $\sum_{k=0}^{+\infty}$ in
{\rref{partpq} is in fact a finite sum $\sum_{k=0}^\ell$.
In the subcase $\ell=0$  of \rref{limizz}, the finite sum consists only of the $k=0$ term, so
\beq {~}_{p} F_{q}(\alpha_1, ... ,\alpha_p; \bet_1,...,\bet_q; \z) =1 \label{limitriv} \feq
$$ \mbox{for $p,q$ arbitrary, if $\alpha_i = 0$ for some $i \in \{1,...,p\}$ and
$\z \in \reali$}~. $$
In general, the series \rref{partpq} is invariant under arbitrary permutations of
the parameters $\alpha_1,...,\alpha_p$ or $\bet_1,...,\bet_q$. \parn
Due to the above indications on the case $p=q$, the function ${~}_{q} F_{q}(\alpha_1,...,\alpha_q;
\bet_1,...,\bet_q; \z)$ is well defined via \rref{partpq} for
\beq \alpha_1,...,\alpha_q \in \reali,~ \bet_1,...,\bet_q \in \reali \setminus (-\naturali),~
w \in \reali~;\label{domanal} \feq
furthermore, ${~}_{q} F_{q}$ is analytic in all the parameters $\alpha_i, \bet_i$ and in the
variable $\z$ on the domain \rref{domanal}.
For fixed $\alpha_1,...,\bet_q$ as in \rref{domanal}, one has
${~}_{q} F_{q}(\alpha_1,...,\alpha_q,\bet_1,...,\bet_q;\z) = O((-\z)^{-\mu})$
for $\z \vain - \infty$, and
${~}_{q} F_{q}(\alpha_1,...,\alpha_q,
\bet_1,...,\bet_q;\z) = O(\z^{\nu} e^\z)$
for $\z \vain + \infty$, where
$\mu := \min(\alpha_1,...,\alpha_q)$,
$\nu := \sum_{i=1}^q \alpha_i - \sum_{i=1}^q \bet_i$;
these results can be traced in the classical work \cite{Bar}. \parn
Concerning the case $p=q+1$, the limitation $\z \in (-1,1)$ in Eq.\rref{limiz} can be
overcome if at least one of the parameters $\alpha_1,...,\alpha_{q+1}$ is positive; in this case,
one can define ${}_{q+1} F_q$ using, instead of the series \rref{partpq},
the following integral formula (see \cite{Luke} Vol.I, page 59, Eq. (13)):
\parn
\vbox{
\beq {~}_{q+1} F_{q}(\alpha_1, ... ,\alpha_{q+1}; \bet_1,...,\bet_q; \z) \label{intef} \feq
$$ :=
{1 \over \Gamma(\alpha_h)} \int_{0}^{+\infty} \!\!\! d t \, e^{-t} \, t^{\alpha_h - 1}
{~}_{q} F_{q}(\alpha_1,...,\alpha_{h-1}, \alpha_{h+1},...\alpha_{q+1}; \bet_1,...,\bet_q; \z t) $$
$$ \mbox{if}~\alpha_h \in (0,+\infty)~\mbox{for some}~h \in \{1,...,q+1\}~\mbox{and}~
\alpha_1,...,\alpha_{h-1}, \alpha_{h+1},...\alpha_{q+1} \in \reali, $$
$$ \bet_1, ..., \bet_q \in \reali \setminus (-\naturali),~ \z \in (-\infty, 1)~. $$}
The above integral converges, due to the previous result on the
asymptotics of ${~}_{q} F_{q}$ for large values of
the variable. The prescription \rref{intef} gives a unique definition for ${~}_{q+1} F_{q}$
if applied for different values of $h$ (all of them with $\alpha_h > 0$), and always
agrees with Eq. \rref{partpq} if $w \in (-1,1)$, or if $\alpha_i = -s$
for some $i \in \{1,...,p\}$, $s \in \naturali$ and $\z \in (-\infty,1)$. \parn
The function ${~}_{q+1} F_{q}$ is analytic in the parameters $\alpha_1,...,\alpha_{q+1},
\bet_1,...,\bet_q$ and in the variable $\z$ in the domain indicated by Eqs.
\rref{limiz} \rref{limizz} and \rref{intef}. Of course, many properties of ${~}_{q+1} F_{q}$
derivable where the series \rref{partpq} converges hold in fact on
the whole domain \rref{limiz} \rref{limizz} \rref{intef}, by the principle of analytic continuation. \parn
Let us finally mention that, for $i \in \{1,...,p\}$ and
$j \in \{1,...,q\}$,
\beq {}_{p+1} F_{q + 1}(\alpha_1,..., \alpha_{i-1}, \gaa, \alpha_{i}..., \alpha_p;
\bet_1,..., \bet_{j-1}, \gaa, \bet_{j}..., \bet_q; \z) \label{ofco} \feq
$$ = {}_{p} F_{q}(\alpha_1,..., \alpha_{i-1}, \alpha_{i}..., \alpha_p;
\bet_1,..., \bet_{j-1}, \bet_{j}..., \bet_q; \z) $$
whenever the two sides are defined (by power series of the type \rref{partpq}, or by
any analytic continuation). \parn
As anticipated, in this paper we are mainly interested in the
${}_2 F_1$ and ${}_{3} F_2$ hypergeometric functions. \parn
The properties of ${}_{2} F_{1}(\alpha,\gammaa;\epp; \z)$
we are using more frequently are the obvious symmetry in $\alpha, \gammaa$,
and the Kummer transformation
\beq {}_{2} F_{1}(\alpha,\gammaa;\epp; \z) = (1 - \z)^{\epp - \alpha - \gammaa}
{}_{2} F_{1}(\epp - \alpha,\epp - \gammaa;\epp; \z)~. \label{sukum} \feq
Besides the integral representation \rref{intef}, we have for this function
the alternative representations \parn
\vbox{
\beq _{2} F_{1}(\alpha, \gammaa; \epp; \x) =
{\Gamma(\epp) \over \Gamma(\gammaa) \Gamma(\epp-\gammaa)}~\int_{0}^{1} d s \, s^{\gammaa-1} (1-s)^{\epp-\gammaa-1}
(1 - s \x)^{-\alpha} \label{irep} \feq
$$ \quad \mbox{for~ $\epp > \gammaa > 0$,~ $-\infty < \x < 1$}~; $$}
\beq {~}_{2} F_{1}(\alpha, \gammaa; \epp; 1-\x) =
{\Gamma(\epp) \over \Gamma(\gammaa) \Gamma(\epp-\gammaa)}~
\int_{0}^{+\infty} \!\!\! du \, u^{\gammaa-1} (1 + u)^{\alpha-\epp} (1 + \x u)^{-\alpha} > 0 \label{gf} \feq
$$ \mbox{for~ $\epp > \gammaa > 0$,~ $\x > 0$}~. $$
Eq. \rref{irep} is the well known Euler's formula, and
\rref{gf} follows from \rref{irep} after a change of variable $s = u/(1+u)$. \parn
The function ${~}_{3} F_{2}(\alpha,\beta,\gamma; \delta, \ep; \w)$
is obviously symmetric in $\alpha, \beta, \gamma$
and $\delta,\ep$ separately.
In the sequel we refer to the identity
(see \cite{Luke}, Vol. II, page 13, Eq. (34))
\parn
\vbox{
$$ {~}_{3} F_{2}(\alpha,\beta,\gamma; \delta, \ep; \z) =
\sum_{i=0}^{+\infty} {(\alpha)_{i} (\beta)_i (\ep - \gamma)_i \over (\delta)_{i} (\ep)_i}
 {(-\z)^{i} \over i!} \, {~}_{2} F_{1}(\alpha + i,\beta + i;\delta+ i; \z) $$
\beq ~\mbox{for $- \infty < \z < {1 \over 2}$~,} \label{dadare} \feq}
We also mention the asymptotics  \cite{KS} \cite{PBM}
\beq {~}_{2} F_{1}(\alpha, \gammaa; \epp; \x) \sim {\Gamma(\gammaa - \alpha) \Gamma(\epp)
\over \Gamma(\epp - \al) \Gamma(\gammaa)} \, (-\x)^{-\al}
\label{asi21} \feq
$$ \mbox{for $\x \vain -\infty$, $~\gammaa, \epp > 0$, $~\al < \min(\gammaa, \epp)$~;} $$
\beq {~}_{3} F_{2}(\alpha,\beta,\gamma; \delta, \ep; \z) \sim
{\Gamma(\de) \Gamma(\ep) \Gamma(\be - \al) \Gamma(\gamma - \al)
\over \Gamma(\beta) \Gamma(\gamma) \Gamma(\de - \al) \Gamma(\ep - \al)}~(-\z)^{-\al}
\label{asimp} \feq
$$ \mbox{for $\z \vain - \infty$, $~~\be,\gamma,\de,\ep > 0$, $~\al < \min(\be,\gamma,\de,\ep)$}~. $$
Another result, important for our purposes, is the relation
\beq \int_{0}^{+\infty} \!\!\! dr\, r^{\mu + \nu + \delta + 1} J_{\delta}(h r) K_{\mu}(r) K_{\nu}(r)  \label{result} \feq
$$ = 2^{\mu + \nu + \delta - 1}~
{\Gamma(\mu + \delta + 1) \Gamma(\nu + \delta + 1) \Gamma(\mu + \nu + \delta + 1) \over
\Gamma(\mu + \nu + 2 \delta + 2 )} \, h^{\delta} $$
$$ \times
{~}_{3} F_{2}(\mu + \delta + 1, \nu + \delta + 1, \mu + \nu + \delta + 1;
{\mu + \nu \over 2} + \delta + 1, {\mu + \nu \over 2} + \delta + {3 \over 2}; -{h^2 \over 4}) $$
$$ \mbox{for $h,\mu,\nu,\delta \in \reali$,~ $h > 0$,
$\delta, \mu + \delta, \nu + \delta, \mu + \nu + \delta > -1$}~; $$
the above conditions on the parameters ensure, amongst else, convergence of the integral
in the left hand side.
Eq.\,\rref{result} generalizes Eq.\,(3.16) of \cite{uno}, and the considerations
of the cited reference can be rephrased in the present framework: the result \rref{result} is
known, but it is difficult to trace a proof in the literature. For this reason,
a derivation of \rref{result} is proposed in Appendix \ref{appef}.
\vskip 0.2cm \noindent
\textbf{Fourier transform.} Let us use the standard notation $S'(\reali^d, \complessi)$
for the tempered distributions on $\reali^d$.
We denote with $\FF, \FF^{-1} : S'(\reali^d, \complessi) \vain S'(\reali^d, \complessi)$ the Fourier transform
and its inverse; $\FF$ is normalized so that
\beq \FF f(k)=
{1 \over (2 \pi)^{d/2}} \int_{\reali^d} d x~e^{-i k \bullet x} f(x) \feq
(intending the integral literally, if $f \in L^1(\reali^d, \complessi)$).
The restriction of $\FF$ to $\Ld$, with the standard inner product and
the associated norm $\|~\|_{L^2}$, is a Hilbertian isomorphism. \parn
Consider two (sufficiently regular) radially symmetric
functions
\beq f : \reali^d \vain \complessi,~~
x \vain f(x) = \varphi(| x |)~, \quad F : \reali^d \vain \complessi,~~
k \vain F(k) = \Phi(| k |)~; \feq
the Fourier and inverse Fourier transforms $\FF f$, $\FF^{-1} F$ are also radially symmetric, and
given by \cite{Boc}
\beq \FF f(k) = {1 \over | k |^{d/2 - 1}}~
\int_{0}^{+\infty} \!\!\! dr~ r^{d/2} J_{d/2 - 1}(| k | r) \varphi(r)~, \label{eboc} \feq
\beq \FF^{-1} F (x) = {1 \over | x |^{d/2 - 1}}~
\int_{0}^{+\infty} \!\!\! d\rho~ \rho^{d/2} J_{d/2 - 1}(| x | \rho) \Phi(\ro)~. \label{ebocc} \feq
\vskip 0.2cm \noindent
\textbf{Sobolev spaces.} Let us consider a real number $\ell$; we denote with $\sqrt{1 + | \k |^2}^{~\ell}$
the function $k \in \reali^d \mapsto \sqrt{1 + | k |^2}^{~\ell}$ (and the multiplication operator
by this function). Furthermore, we put
\beq \sqrt{1 - \Delta}^{~\ell}~ := \FF^{-1} \sqrt{1 + | \k |^2}^{~\ell} \FF :  S'(\reali^d, \complessi)
\vain S'(\reali^d, \complessi)~. \label{lap} \feq
The $\ell$-th order Sobolev (or Bessel potential) space of
$L^2$-type and its norm are \cite{Smi} \cite{Maz}
\beq H^\ell(\reali^d, \complessi) := \lgraffa f \in
S'(\reali^d, \complessi)~\Big\vert~ \sqrt{1 - \Delta}^{~\ell} f \in
\Ld~ \rgraffa
\label{incid}\feq
$$ = \lgraffa f \in S'(\reali^d, \complessi)~\Big \vert~
\sqrt{1 + | \k |^2}^{~\ell} \FF f \in \Ld \rgraffa~;
$$
\beq \| f \|_{\ell} := \| \sqrt{1 - \Delta}^{~\ell}~ f \|_{L^2} =
\|\, \sqrt{1 + | \k |^2}^{~\ell}~\FF f \, \|_{L^2}~. \label{repfur} \feq
We note the equality
\beq (H^0(\reali^d), \|~\|_0) = (L^2(\reali^d), \|~\|_{L^2}) \feq
and the imbedding relations
\beq \ell \leqs \ell' \qquad \Rightarrow \qquad H^{\ell'}(\reali)^d \subset H^\ell(\reali^d)~,\quad
\| ~\|_\ell \leqs \| ~\|_{\ell'}~. \feq
We only consider the Sobolev spaces $H^\ell(\reali^d)$
of order $\ell \geqs 0$, which are embedded into $L^2(\reali^d)$ (and so, consist of ordinary
functions).
In the special case $\ell \in \naturali$, the definitions \rref{incid} \rref{repfur} imply
\parn
\vbox{
\beq H^\ell(\reali^d, \complessi) = \{ f \in S'(\reali^d, \complessi)~|~
\partial_{\lambda_1,..., \lambda_k} f \in L^2(\reali^d, \complessi)~\label{hnab} \feq
$$ \forall k \in \{0, ..., \ell\}, (\lambda_1, ..., \lambda_k) \in \{1, ..., d\}^{k}~\}~; $$
\beq \| f \|_\ell = \sqrt{\sum_{k=0}^\ell \left( \barray{c} \ell \\ k \farray \right)
\sum_{\lambda_1, ... ,\lambda_k =1, ..., d}
\int_{\reali^d} d x~| \partial_{\lambda_1,..., \lambda_k} f (x) |^2}~. \label{nonab} \feq}
In the above, $\partial_{\lambda_i}$ is the distributional derivative with respect to the coordinate
$x^{\lambda_i}$.
\vskip 0.2cm \noindent
\textbf{Other functions.} As in \cite{uno}, a central role in our considerations
is played by the function $G_{t d} := 1/(1 + |\k|^2)^{t}$, i.e.,
\beq \F_{t d} : \reali^d \vain \complessi~,
\quad k \mapsto \F_{t d}(k) := {1 \over (1 + | k |^2)^t} \qquad (t \in \reali)~; \label{efnd} \feq
we further set
\beq \f_{t d} : \reali^d \vain \complessi~, \qquad \f_{t d} := \FF^{-1} \F_{t d}
\qquad (t > d/4)~.  \label{fnd} \feq
We note that, with the assumption $t > d/4$, $\F_{t d}$ and, consequently, $\f_{t d}$ are
$L^2$ functions.
The functions $\f_{t d}$ are related to the Macdonald functions \cite{Smi} \cite{Maz} since, for any
$x \in \reali^d$,
\beq \f_{t d}(x) = {| x |^{t - d/2} \over 2^{t - 1} \Gamma(t)}~
K_{t - d/2}(| x |)~. \label{gemac} \feq
\section{The constant $\boma{K_{\ell m n d}}$ and
its bounds: description of the main results.}
\label{desc}
Let $d \in \naturali_0$, and consider three real numbers $\ell, m, n$ such that
\beq 0 \leqs \ell \leqs m \leqs n~,\qquad n + m - \ell > d/2~. \label{c1} \feq
\begin{prop}
\label{demul}
\textbf{Definition.} We put
\beq K_{\ell m n d} := \min~\Big\{~ K \in [0,+\infty)~\Big|~ \| f g \|_\ell \leqs K \| f \|_m \| g \|_n  ~~\label{bemul} \feq
$$\mbox{for all $f \in \Hm, g \in \Hn$}~\Big \}$$
and refer to this as the sharp (or best) constant for the multiplication $\Hm \times \Hn \vain \Hell$.
\end{prop}
In the sequel we present our upper and lower bounds for the above constant; most of
the forthcoming propositions are proved in Sections \ref{upper}, \ref{pbes}, \ref{secdir}.
\vskip 0.2cm \noindent
\textbf{``$\boma{\SS}$-function'' upper bound on} $\boma{K_{\ell m n d}}$. This is our most important
upper bound; it is determined by a function $\SS = \SS_{\ell m n d}$, as stated hereafter.
\begin{prop}
\label{pupper}
\textbf{Proposition.} (i) For $\ell, m, n$ fulfilling \rref{c1}, one has
\beq K_{\ell m n d} \leqs \sqrt{\sup_{u \in [0,+\infty)} \SS_{\ell m n d}(u)}~,
\label{kpnd} \feq
where, for $u \in [0,+\infty)$,
\beq \SS_{\ell m n d}(u) := {\Gamma(m + n - d/2) \over (4 \pi)^{d/2} \Gamma(n + m)}
(1 + 4 u)^\ell~F_{m n d}(u)~, \label{ff} \feq
\beq F_{m n d}(u) := {~}_{3} F_{2}(m + n - {d \over 2}, m, n; {m + n \over 2}, {m + n + 1 \over 2}; -u)~.
\label{fmnd} \feq
In the special case $m=n$, Eq. \rref{fmnd} implies
\beq F_{m m d}(u) = {~}_{2} F_{1}(2 m - {d \over 2}, m; m + {1 \over 2}; -u)~;
\label{fmmd} \feq
the trivial case $m=0$ is described by
\beq F_{0 n d}(u) = 1 ~\qquad \mbox{for all $u$}~. \label{ftriv} \feq
For all $\ell, m, n$ as in \rref{c1},
the function $\SS_{\ell m n d}$ sends $[0,+\infty)$ to $(0,+\infty)$ and is bounded, so
the sup in \rref{kpnd} is actually finite. The behavior of this
function for $u=0$ and $u \vain + \infty$ is described by the following relations:
\beq \SS_{\ell m n d}(0) = {\Gamma(m + n - d/2) \over (4 \pi)^{d/2} \Gamma(n + m)}~, \label{duesei} \feq
\beq \SS_{\ell m n d}(u) \sim {(1 + \delta_{m n}) \Gamma(n - d/2) \over (4 \pi)^{d/2} \Gamma(n)}
~{1 \over (4 u)^{m - \ell}} \qquad \mbox{for $u \vain + \infty$} \label{dueset} \feq
($\delta$ is the Kronecker symbol, i.e., $\delta_{m n} := 1$ if $m=n$, and
$\delta_{m n} :=0$ otherwise).
According to \rref{dueset}, the $u \vain +\infty$ limit of $\SS_{\ell m n d}$ is
\beq \SS_{\ell m n d}(+\infty) = \left\{ \barray{ll} \dd{(1 + \delta_{m n}) \Gamma(n - d/2) \over (4 \pi)^{d/2} \Gamma(n)}
& \mbox{if $\ell = m$,} \\ 0 & \mbox{if $\ell < m$.}\farray \right.  \label{dueot} \feq
\vbox{\noindent
(ii) One has
\beq F_{m n d}(u) \label{esp1} \feq
$$ = \sum_{i=0}^{+\infty} \sum_{j=0}^{+\infty} { \left(m + n - {d \over 2} \right)_{i} (m)_i
\left( {m - n + 1 \over 2} \right)_{i} \over i! \left({m + n + 1 \over 2} \right)_{i}
\left({m + n  \over 2} \right)_{i} }~{\left( {d - m - n \over 2} \right)_j
\left( {n - m \over 2} \right)_j \over j! \left( {m + n \over 2} + i \right)_j }
{(-1)^j u^{i+j} \over (1 + u)^{{3 m + n - d \over 2} + i}}~ $$
if $u \in [0,1)$, or $u \in [0,+\infty)$ and the series over $j$ is a finite sum.}
An alternative expansion, holding under the same conditions, is
\parn
\vbox{
\beq F_{m n d}(u) \label{esp2} \feq
$$ = \sum_{i=0}^{+\infty} \sum_{j=0}^{+\infty} { \left(m + n - {d \over 2} \right)_{i} (m)_i
\left( {m - n \over 2} \right)_{i} \over i! \left({m + n + 1 \over 2}
\right)_{i} \left({m + n  \over 2} \right)_{i} }~
{\left( {d + 1 - m - n \over 2} \right)_j
\left( {n + 1 - m \over 2} \right)_j \over j! \left( {m + n + 1 \over 2} + i \right)_j }
{ (-1)^j u^{i+j} \over (1 + u)^{{3 m + n - d - 1\over 2} + i}} ~. $$}
The above series over $j$ or $i$ become finite sums
in the special cases indicated below.
\beq \mbox{If}~~m + n - d \in 2 \naturali, \quad \sum_{j=0}^{+\infty} \vain
\! \sum_{j=0}^{{m + n - d \over 2}}~\mbox{in \rref{esp1}}~;
\label{spe1} \feq
$$ \mbox{if}~~n - m \in 2 \naturali + 1, \quad
\! \sum_{i=0}^{+\infty} \vain \sum_{i=0}^{{n - m - 1 \over 2}} \mbox{in \rref{esp1}}~. $$
\beq \mbox{If}~~m + n - d \in 2 \naturali + 1, \quad \sum_{j=0}^{+\infty} \vain
\! \sum_{j=0}^{{m + n - d - 1 \over 2}} \mbox{in \rref{esp2}}~; \label{spe2} \feq
$$ \mbox{if}~~n - m \in 2 \naturali, \quad
\sum_{i=0}^{+\infty} \vain \sum_{i=0}^{{n - m \over 2}}~~\mbox{in \rref{esp2}}~. $$
\end{prop}
\textbf{Proof.} See Section \ref{upper}. \fine
\begin{rema}
\textbf{Remark.} In the case $\ell=m=n$ ($\ell > d/2$), Eqs. (\ref{ff}-\ref{fmmd}) give
\beq \SS_{\ell \ell \ell d}(u) := {\Gamma(2 \ell - d/2) \over (4 \pi)^{d/2} \Gamma(2 \ell)}
(1 + 4 u)^\ell~{}_{2} F_{1}(2 \ell - {d \over 2}, \ell; \ell + {1 \over 2}; -u)~; \label{ffel} \feq
this is the function denoted with $\SS_{\ell d}$ in \cite{uno}, Proposition 2.2, that was employed
to derive our upper bound on $K_{\ell \ell \ell d} \equiv K_{\ell d}$. \fine
\end{rema}
\salto
\textbf{``H\"older'' upper bound on $\boma{K_{0 m n d}}$.}
The upper bound on $K_{\ell m n d}$ given by the above proposition holds for arbitrary $\ell, m, n$ as in
\rref{c1}. In this paragraph we give a different
upper bound for the special case $\ell=0$, that is somehow trivial since $\|~\|_{0}$ is the
$L^2$-norm. In this case, for all functions $f, g$ one can estimate $\| f g \|_{L^2}$ via the H\"older inequality,
and then employ the Sobolev imbedding inequality, with certain information on the related constant.
To make contact with the Sobolev imbedding, we introduce
the following notations:
\beq R_{\n d} := \left\{\barray{ll} [2, \dd{d \over d/2 - \n} ] &  \mbox{if $\n \in [0, d/2)$}, \\
\mbox{$[2,+\infty)$} &  \mbox{if $\n = d/2$}~,
\\ \mbox{$[2,+\infty]$} & \mbox{if $\n \in (d/2, + \infty)$~;} \farray \right.~ \label{pert} \feq
\beq S_{r \n d} :=
{1 \over (4 \pi)^{d /4 - d/(2 r)} }~
\left( { \Gamma \left( {\displaystyle{ {\n \over 1 - 2/r} - {d \over 2} }} \right)
\over \Gamma
\Big( {\displaystyle{ {\n \over 1 - 2/r} }} \Big) }
\right)^{1/2 - 1/r} \!\!\!\!\!\!\!\! \left({E(1/r) \over E(1 - 1/r)} \right)^{d/2} \label{srid} \feq
$$ \mbox{if~ $\n \in [0, d/2)$, $r \in \big(2, \dd{d \over d/2 - \n}\big)$ ~~or~~ $\n \in [d/2, +\infty)$,
$r \in (2,+ \infty)$}~,
$$
\beq S_{2 \n d} :=  1 \qquad \mbox{if $\n \in [0,+\infty)$}~, \label{s2nd} \feq
\beq S_{\infty \n d} :=
{1 \over (4 \pi)^{d /4} }~
\sqrt{ \Gamma( \n  - d/2 ) \over \Gamma(\n) }
~~\mbox{if $\n \in (\dd{d / 2}, + \infty)$}~; \label{siid} \feq
\beq E(u) := u^u \qquad \mbox{for $u \in (0, + \infty)$}~,
\qquad E(0) := \mbox{lim}_{u \vain 0^{+}} E(u) = 1~. \label{siie} \feq
Then
\beq \n \in [0,+\infty], r \in R_{\n d} \quad\Rightarrow\quad H^\n(\reali^d) \subset L^r(\reali^d),
~\|~\|_{L^r} \leqs S_{r \n d} \, \|~ \|_{\n}~; \label{imbrel} \feq
furthermore, for $\n \in (d/2, + \infty)$,
\beq S_{\infty \n d} := \min \{ S \in [0,+\infty)~|~~
\|~ \|_{L^\infty(\reali^d)} \leqs S \, \| ~\|_\n~\}. \label{imbinf} \feq
Of course, the imbedding inequality $\|~\|_{L^r(\reali^d)}
\leqs$ constant $\| ~\|_\n$ is well
known; for the statements (\ref{pert}-\ref{imbinf}) on the constant in this inequality,
see \cite{imb}. In particular, \rref{imbinf} means that $S_{\infty \n d }$ is the
sharp constant for the corresponding inequality; as a matter of fact, the equality
$\| f \|_{L^\infty(\reali^d)} = S_{\infty \n d } \, \| f \|_\n$ holds for
$f = \f_{\n d}$ as in Eqs. \rref{fnd} \rref{gemac}.
 \parn
With the above notations, we can state the following.
\parn
\vbox{
\begin{prop}
\label{uppers}
\textbf{Proposition.} For any $p \in [2,+\infty]$, let $p^{*} \in [2,+\infty]$ denote the solution
of the equation $1/p + 1/p^{*} = 1/2$. Furthermore,
let $m, n$ fulfill conditions \rref{c1}, with $\ell=0$; then, (i)(ii)hold. \parn
(i) The set
\beq R_{m n d} := \{ p \in R_{m d}~|~p^{*} \in R_{n d} \} \label{kimbsu0} \feq
is nonempty. \parn
(ii) For any $p \in R_{m n d}$, one has
\beq K_{0 m n d} \leqs S_{p m d} \, S_{p^* n d}~; \label{kimb} \feq
so,
\beq K_{0 m n d} \leqs \inf_{p \in R_{m n d}}
S_{p m d} \, S_{p^{*} n d}~.
\label{kimbsu} \feq
\end{prop}
}
\noindent
\textbf{Proof.}
(i) The thesis follows from an elementary analysis, explicitating
the definitions of $R_{m d}$ and $R_{n d}$ via Eq. \rref{pert}. \parn
(ii) Let $p \in R_{m n d}$, and consider any two functions
$f \in H^m(\reali^d)$, $g \in H^n(\reali^d)$; then, the H\"older inequality and the imbedding
relations \rref{imbrel} give
\beq \| f g \|_{0} = \| f g \|_{L^2} \leqs \| f \|_{L^p} \| g \|_{L^{p^{*}}} \leqs (S_{p m d} \| f \|_m)
(S_{p^{*} n d} \| g \|_n)~, \feq
whence the thesis \rref{kimb}. Now, \rref{kimbsu} is obvious. \fine
As shown later via a series of examples, the bound \rref{kimbsu} is often better than the case $\ell=0$ of the bound
\rref{kpnd}.
\vskip 0.2cm \noindent
\textbf{General method to get lower bounds on} $\boma{K_{\ell m n d}.}$
The general method is based on the obvious inequality
\beq K_{\ell m n d} \geqs {\| f g \|_\ell \over \| f \|_m \| g \|_n} \label{ofc} \feq
for all nonzero $f \in \Hm$, $g \in \Hn$; this gives a lower bound for any pair of ``trial functions'' $f, g$.
In the sequel we propose several choices of the trial functions, depending on one or more parameters;
the parameters must be tuned to get the best lower bound, i.e., the maximum value for the
right hand side of Eq. \rref{ofc}. \parn
\vskip 0.2cm \noindent
\textbf{``Bessel'' lower bound.} In this approach, the trial functions have the form
\beq \f_{\nu t d}(x) := \f_{t d}(\nu x) \label{lagemac} \feq
where $\nu \in (0, +\infty)$ is a parameter and $\f_{t d}$ is defined by Eq.
\rref{fnd}. By comparison with that equation, we find
\beq \f_{\nu t d} = \FF^{-1} \F_{\nu t d}~, \qquad \F_{\nu t d}(k) :=
{1 \over \nu^d (1 + | k |^2/\nu^2)^t}~. \feq
\begin{prop}
\label{pbessel}
\textbf{Proposition.} (i) Let
$n \in [0,+\infty), t \in (\dd{n / 2} + \dd{d / 4}, +\infty)$,
$\nu \in (0,+\infty)$. Then
\beq \f_{\nu t d} \in \Hn,~~\label{fnutd} \feq
$$ \| \f_{\nu t d} \|^2_{n} = {\pi^{d / 2} \over \nu^d}
{\Gamma(2 t - n - d/2) \over \Gamma(2 t - n)} {~}_2 F_1(- n, {d / 2}; 2 t - n; 1 - \nu^2)~. $$
(Note that ${~}_2 F_1(- n, \dd{d / 2}; 2 t - n; w)$ is a finite sum
$\sum_{i=0}^{n} \dd{ (-n)_i (d / 2)_i \over (2 t - n)_i}
\dd{w^i \over i!}$ if $n \in \naturali$). \parn
(ii) Let $\ell, m, n$ fulfill \rref{c1}, and
\beq s \in (\dd{m / 2} + \dd{d / 4}, +\infty)~, \quad t \in (\dd{n / 2} + \dd{d / 4}, +\infty)~,
\qquad \mu, \nu \in (0,+\infty)~ \label{dacit} \feq
(then $\f_{\mu s d} \in H^m(\reali^d, \complessi)$
and $\f_{\nu t d} \in H^n(\reali^d, \complessi)$, due to (i); this also implies
$\f_{\mu s d} \, \f_{\nu t d} \in \Hell$). One has
\beq \| \f_{\mu s d} \, \f_{\nu t d} \|^2_{\ell} = {2^d \pi^{d/2} \over \Gamma(d/2)}
\int_{0}^{+\infty} \!\!\! du \,
u^{d/2 - 1} (1 + 4 u)^{\ell} G^2_{s t d}(\mu, \nu; u)~, \label{fmufnu} \feq
where
\beq G_{s t d}( \mu, \nu; u) \label{gstd} \feq
$$ := {\mu^{s-d/2} \nu^{t-d/2} \over {2^{2 s + 2 t - 2} \Gamma(s) \Gamma(t) u^{s/2 + t/2}}}~
\int_{0}^{+\infty} \!\!\! dr \, r^{s + t - d/2} J_{d/2 - 1}(r)~K_{s - d/2}({\mu r \over 2 \sqrt{u}})
K_{t - d/2}({\nu r \over 2 \sqrt{u}}) ~. $$
Moreover, assume
\beq s - {d \over 2}, \, t - {d \over 2} \in \naturali + {1 \over 2}, ~~\ell \in \naturali. \label{condst} \feq
Then both integrals in Eqs. \rref{gstd} and \rref{fmufnu} are elementary, and
\beq \| \f_{\mu s d} \, \f_{\nu t d} \|^2_{\ell} =
{\pi^{d/2 + 2} \over \Gamma^3(d/2) \Gamma^2(s) \Gamma^2(t)} \,
\sum_{h=0}^{\ell} \sum_{(i,j,k) \in I_{s t d}} \sum_{(i',j',k') \in I_{s t d}}
\left( \barray{cc} \ell \\ h \farray \right)
\label{pooint} \feq
$$ \times {
\Gamma (i + i' + j+ j' - k - k' - h +  d/2 + 1) \Gamma (k + k' + h + d/2)
\over \Gamma(i+i' + j+j' + d + 1)}\,
G_{s t i j k d} \, G_{s t i' j' k' d} $$
$$ \times {\mu^{i+i'} \nu^{j+j'} \over (\mu + \nu)^{i + i' + j + j' - 2 h +  d}}~. $$
Here we have put
\beq I_{s t d} \label{istd} \feq
$$ := \{ (i, j, k) \in \naturali^3~|~0 \leqs i \leqs s-{d \over 2} - {1 \over 2},
0 \leqs j \leqs t- {d \over 2} - {1 \over 2}, 0 \leqs k \leqs {i + j + 1 \over 2}~\}~; $$
\noindent
\vbox{
\beq G_{s t i j k d} \label{gstijkd} \feq
$$ :=
{(-1)^k (i + j + d-1)! (2 s - i - d - 1)! (2 t - j - d - 1)!
\left(-{i + j \over 2}\right)_k \left(-{i + j + 1 \over 2}\right)_k
\over 2^{2 s + 2 t - i - j - d/2 - 3} \, i! \, j! \, k!
(s - i - {d \over 2} - {1 \over 2})! \, (t - j - {d \over 2} - {1 \over 2})!
\left({d \over 2}\right)_k}~. $$
}
\noindent
(iii) Let $\ell, m, n$ be as in \rref{c1}, and $s,t$ as in (ii). Then, for all $\mu, \nu \in(0,+\infty)$,
\beq K_{\ell m n d} \geqs \Kap^{\Bst}_{\ell m n d}(\mu, \nu) :=
{ \| \f_{\mu s d} \f_{\nu t d} \|_\ell \over \| \f_{\mu s d} \|_m \| \f_{\nu t d}\|_n}~,
\label{theabove} \feq
whence
\beq K_{\ell m n d} \geqs \sup_{\mu,\nu > 0}
\Kap^{\Bst}_{\ell m n d}(\mu, \nu)~. \label{thea1} \feq
The function $\Kap^{\Bst}_{\ell m n d}$ can be computed from items (i)(ii).
\end{prop}
\parn
\textbf{Proof.} See Section \ref{pbes}. \fine
\vskip 0.2cm \noindent
\textbf{``Fourier" lower bound on $\boma{K_{\ell m n d}}$}.
As in \cite{uno}, we use this term for the lower bound arising from the trial functions
\beq f_{p \sigma d}(x) := e^{i p x_1}~e^{- \sigma | x |^2/2}
\qquad (p \in [0,+\infty), \sigma \in (0,+\infty)) \label{fcar} \feq
\parn
The Sobolev norm of any order $n$ of this
function can be expressed using the modified Bessel function of the first kind
$I_{\nu}$, the Pochhammer symbol \rref{poch} and the double factorial \rref{semif}.
\vskip 0.2cm \noindent
\begin{prop}
\label{pfou} \textbf{Proposition.}
(i) Let $m, p \in [0,+\infty)$, $\sigma \in (0,+\infty)$. Then
\beq \| f_{p \sigma d} \|_m^2 = {2 \, \pi^{d/2}  \over \sigma^{d/2 + 1} p^{d/2 - 1}}~
\int_{0}^{+\infty} \!\!\! d\rho \,\rho^{d/2} (1 + \rho^2)^m e^{-{\rho^2 + p^2\over \sigma}} I_{d/2 - 1} ({2 p \over
\sigma} \rho) \label{giveby} \feq
if $p > 0$, and
\beq \| f_{0 \sigma d} \|_m^2 = {2 \, \pi^{d/2}  \over \Gamma(d/2) \sigma^{d}}~
\int_{0}^{+\infty} \!\!\! d\rho \,\rho^{d-1} (1 + \rho^2)^m e^{-{\rho^2 \over \sigma}}~ \label{giveby0} \feq
(this is the $p \vain 0^{+}$ limit of \rref{giveby}, since $I_{d/2 -1}(w) \sim \dd{(w/2)^{d/2-1} \over \Gamma(d/2)}$
for $w \vain 0^{+}$). \parn
In particular, for $m$ integer,
$$ \| f_{p \sigma d} \|_m^2 = \pi^{d/2} \sum_{\ell=0}^m \sum_{j=0}^\ell \sum_{g=0}^j
\left( \barray{c} m \\ \ell \farray \right) \left( \barray{c} \ell \\ j \farray \right)
\left( \barray{c} 2 j \\ 2 g \farray \right) {(2 g - 1)!! \over 2^g}  $$
\beq \times \left({d/2}- {1/2}\right)_{\ell-j}
p^{2 j - 2 g} \sigma^{\ell + g - j - d/2}~. \label{gv} \feq
(ii) Let $\ell, m, n$ fulfill \rref{c1}. Then, for all $p, q \in [0,+\infty)$
and $\si, \tau \in (0,+\infty)$,
\beq K_{\ell m n d} \geqs \Kap^{F}_{\ell m n d}(p, q, \si, \tau) :=
{ \| f_{p + q, \si + \tau, d} \|_\ell \over \| f_{p \si  d} \|_m \| f_{q \tau d}\|_n}~,
\label{theabovee} \feq
whence
\beq K_{\ell m n d} \geqs \sup_{p, q \geqs 0, \, \si, \tau > 0}
\Kap^{F}_{\ell m n d}(p, q, \si, \tau)~. \label{thea11} \feq
The function $\Kap^{F}_{\ell m n d}$ can be computed from item (i).
\end{prop}
\textbf{Proof.} (i) See \cite{uno}, Proposition 2.4. \parn
(ii) Use Eq. \rref{ofc} with $f = f_{p \si d}$ and $g = f_{q \tau d}$; then
$f g = f_{p + q, \si + \tau, d}$ and we get Eq. \rref{theabovee}.
\fine
\vskip 0.2cm \noindent
\textbf{``$\boma{S}$-constant'' lower bound on $\boma{K_{\ell \ell n d}}$.}
This lower bound holds for $K_{\ell m n d}$ in the special case $\ell =m$;
it can be obtained from \rref{ofc}, substituting for $f$ a family of approximants of the
Dirac $\delta$ distribution. This bound already appeared in \cite{mp2}, analyzing an inequality
strictly related to the case $\ell=m$ of \rref{bemul}. In the cited reference,
for a number of reasons this was called  the ``ground level'' lower bound; here,
we prefer the denomination of ``$S$-constant'' lower bound to recall its
relation with the Sobolev imbedding constant $S= S_{\infty n d}$ of Eq. \rref{siid}.
\begin{prop}
\label{propdir}
\textbf{Proposition.} Let
\beq 0 \leqs \ell \leqs n~, \qquad n > {d \over 2}~. \label{lett} \feq
Then
\beq K_{\ell \ell n d} \geqs S_{\infty n d}~. \label{kell} \feq
\end{prop}
\textbf{Proof.} It is essentially known from \cite{mp2}; for completeness, a
sketch of it is given in Section \ref{secdir}. \fine
The last statement, combined with the general upper bound \rref{kpnd} in
Proposition \ref{pupper}, gives
the sharp value of $K_{\ell \ell n d}$ in the trivial case $\ell=0$.
\begin{prop}
\textbf{Proposition.} Let $n > d/2$; then
\beq K_{0 0 n d} = S_{\infty n d}~. \label{genthe} \feq
\end{prop}
\textbf{Proof.}
The cited inequality \rref{kpnd} gives
\beq K_{0 0 n d} \leqs \sqrt{\sup_{u \in [0,+\infty)} \SS_{0 0 n d}(u)}~; \label{kpnd0} \feq
on the other hand, the general definition \rref{ff} of $\SS_{\ell m n d}$
and Eq. \rref{ftriv} about $F_{0 n d}$ give
\beq \SS_{0 0 n d}(u) = {\Gamma(n - d/2) \over (4 \pi)^{d/2} \Gamma(n)} F_{0 n d}(u)
= {\Gamma(n - d/2) \over (4 \pi)^{d/2} \Gamma(n)}
\qquad \mbox{for all $u \in [0,+\infty)$}~. \feq
From here, \rref{kpnd0} and \rref{siid} we see that
\beq K_{0 0 n d} \leqs {1 \over (4 \pi)^{d/4}} \sqrt{\Gamma(n - d/2) \over \Gamma(n)} =
S_{\infty n d}~. \feq
From \rref{kell} we have $K_{0 0 n d} \geqs S_{\infty n d}$ as well, so we get
the thesis \rref{genthe}. \fine
In fact, the equality $K_{\ell \ell n d} = S_{\infty n d}$ holds as well in some cases
with nonzero $\ell$ (e.g., for $d=3$ and $\ell =1, n=2$: see the table of page
\pageref{patable} and Eqs. (\ref{eqqa}-\ref{eqqc})).
\vskip 0.2cm \noindent
\section{On the explicit determination of upper and lower bounds for
$\boma{K_{\ell m n d}}$\,.}
\label{sectable}
Let us translate the results of the previous section into a scheme
to get explicit upper and lower bounds $K^{\pm}_{\ell m n d}$
on $K_{\ell m n d}$, such that
$$ K^{-}_{\ell m n d} \leqs K_{\ell m n d} \leqs K^{+}_{\ell m n d}~. $$
At the end of the section, we present a table of such upper and lower bounds,
for $d=1$ or $3$ and many values of $\ell, m, n$. Before discussing the table, let us describe
the general scheme to determine the upper and lower bounds. \parn
\vskip 0.2cm \noindent
\vbox{
\noindent
\textbf{On the computation of $\boma{K^{+}_{\ell m n d}}$}. One proceeds as follows. \parn
(i) For any $\ell \geqs 0$, one can use the $\SS$-function upper bound provided by Proposition \ref{pupper}, Eq.
\rref{kpnd}, i.e., the number}
\beq  K^{\SS}_{\ell m n d} :=
\sqrt{\sup_{u \in [0,+\infty)} \SS_{\ell m n d}(u)}~~\mbox{(or an upper approximant for this)}. \label{kpnnd} \feq
The function $\SS_{\ell m n d}$
has the expression provided by Eqs. (\ref{ff}-\ref{spe2}); depending on the case,
its $\sup$ can be determined analytically or estimated numerically. \parn
(ii) For $\ell = 0$, one can use as well
the H\"older upper bound provided by Proposition \ref{uppers}, Eq. \rref{kimb}, i.e., the number
\beq  K^{\Ii}_{0 m n d} := \inf_{p \in R_{m n d}}
S_{p m d} \, S_{p^* n d}~~\mbox{(or an upper approximant for this)}~. \label{kimbsuu} \feq
Let us recall that $1/p + 1/p^{*} =1/2$ and $S_{p m d}$, $R_{p m d}$ are
defined by Eqs. (\ref{srid}-\ref{siie}), \rref{kimbsu0}; typically,
the estimation of the sup over $p$ is numerical. \parn
(iii) We denote with $K^{+}_{\ell m n d}$ the best upper bound arising from (i) (ii); so
\beq K^{+}_{\ell m n d} := K^{\SS}_{\ell m n d}~~\mbox{if $\ell > 0$}~, \qquad
K^{+}_{0 m n d} := \min(K^{\SS}_{0 m n d}, K^{\Ii}_{0 m n d})~. \label{minkp} \feq
\vskip 0.2cm \noindent
\textbf{On the computation of $\boma{K^{-}_{\ell m n d}}$}\,.
One proceeds in this way (possibly using numerical methods to compute the
quantities mentioned below). \parn
(i) One chooses two values $(s,t)$ fulfilling conditions \rref{dacit}; the choice
$s = m, t = n$ is natural whenever possible. After fixing $s, t$ one considers
for $K_{\ell m n d}$ the Bessel lower bound suggested by Proposition \ref{pbessel}, Eq.
\rref{thea1}, i.e., the number
\beq  K^{\Bst}_{\ell m n d} := \sup_{\mu, \nu > 0} \Kap^{\Bst}_{\ell m n d}(\mu, \nu)
~~\mbox{(or a lower approximant for this)}~. \label{kbes} \feq
The function $\Kap^{\Bst}_{\ell m n d}$ is
determined by Eqs. (\ref{fnutd}-\ref{theabove}). \parn
(ii) An alternative to the bound \rref{kbes} is the
Fourier lower bound suggested by Proposition \ref{pfou}, Eq. \rref{thea11}, i.e.,
the number
\beq K^{F}_{\ell m n d} := \sup_{p,q \geqs 0, \, \si,\tau > 0}
\Kap^{F}_{\ell m n d}(p, q, \si, \tau)
~~\mbox{(or a lower approximant for this)}~. \label{thea111} \feq
The function $\Kap^{F}_{\ell m n d}$ is determined by
Eqs. (\ref{giveby}-\ref{theabovee}). \parn
(iii) In the special case $\ell=m$, Proposition \ref{propdir} also gives the $S$-constant lower bound
$$ K_{\ell \ell n d} \geqs S_{\infty n d}~, $$
with $S_{\infty n d}$ as in \rref{siid}. \parn
(iv) The best lower bound arising from (i) (ii) (iii) is
\beq K^{-}_{\ell m n d} := \max(K^{\Bst}_{\ell m n d}, K^{F}_{\ell m n d})~
\mbox{if $\ell < m$},~~
K^{-}_{\ell \ell n d} := \max(K^{\Bst}_{\ell \ell n d}, K^{F}_{\ell \ell n d}, S_{\infty n d}). \label{kmaxx} \feq
\textbf{A table of upper and lower bounds.} The forthcoming table considers the dimensions
$d=1,3$ and a set of integer values for $\ell, m, n$. For each one of these values
an upper bound $K^{+}_{\ell m n d}$ and a lower bound $K^{-}_{\ell m n d}$ have been
computed with the methods outlined above. Then, the values of $K^{+}_{\ell m n d}$ and of the ratio
$K^{-}_{\ell m n d}/K^{+}_{\ell m n d}$ have been reported in the table: giving the above ratio,
rather than the lower bound, is more convenient to appreciate
how narrow is the uncertainty on $K_{\ell m n d}$. \parn
In all cases considered in the table
$\SS_{\ell m n d}$, $\Kap^{\Bst}_{\ell m n d}$ and $\Kap^{F}_{\ell m n d}$ are elementary
functions, but often they have lengthy expressions; typically,
their sups or infs have been evaluated numerically. The long expressions for the cited
functions have been obtained implementing the general formulas
of Section \ref{desc} on MATHEMATICA, in the symbolic mode; the same package, with its standard optimization
algorithms, has been employed to compute numerically the necessary sups and infs. \parn
In the cases $\ell=0$ of the table, the minimum \rref{minkp} giving
$K^{+}_{0 m n d}$ equals $K^{\Ii}_{0 m n d}$. \parn
Depending on the case, the
lower bound $K^{-}_{\ell m n d}$ in \rref{kmaxx} can either be
a Bessel bound $K^{\Bst}_{\ell m n d}$, a Fourier bound
$K^{F}_{\ell m n d}$ or an $S$-constant bound $S_{\infty n d}$; to distinguish these situations we
have placed after the value
of $K^{-}_{\ell m n d}/K^{+}_{\ell m n d}$ the symbols $(\Bst)$, $(F)$ or $(S)$, respectively. \parn
Hereafter we present the table of upper and lower bounds; in Appendix \ref{apptab}
we give some examples of the calculations from which the table originated,
reporting all the necessary details.
\vfill \eject \noindent
\oddsidemargin=0.3truecm
\vskip -1cm \noindent
\textbf{Table of the bounds $\boma{K^{-}_{\ell m n d} \leqs K_{\ell m n d} \leqs K^{+}_{\ell m n d}}$
for $\boma{d=1,3}$ and
some values of $\boma{\ell, m, n}$}
\textbf{(the notations $\boma{(F), (B_{s t}), (S)}$ indicate the type of the lower bound
$\boma{K^{-}_{\ell m n d}}$).} \label{patable}
\vskip 0.3cm \noindent
\hrule
\vskip 0.4cm \noindent
$\boxed{\boma{d=1}}$ \hskip 6.3cm $\boxed{\boma{d=3}}$
\vskip 0.1cm \noindent
\renewcommand\tabcolsep{0.1cm}
\renewcommand\arraystretch{1.3}
{\footnotesize{
$$ \begin{tabular}{c|c|c||c|cl}
$\ell$ & $m$ & $n$ & ~~$K^{+}_{\ell m n d}~~ $ & $K^{-}_{\ell m n d} /K^{+}_{\ell m n d}$ \\
\hline \hline
$0$ & $1$ & $1$ & 0.439 & \hs 0.917 & \hs ($B_{1 1}$) \\
\hline
$0$ & $1$ & $2$ & 0.383 & \hs 0.987 & \hs ($F$) \\
\hline
$0$ & $1$ & $10$ & 0.274 & \hs 0.997 & \hs ($F$) \\
\hline
$1$ & $1$ & $2$ & 0.562 & \hs 0.916 & \hs ($B_{1 2}$) \\
\hline
$1$ & $1$ & $3$ & 0.464 & \hs 0.945 & \hs ($B_{1 3}$) \\
\hline
$1$ & $1$ & $10$ & 0.310 & \hs 0.984 & \hs ($B_{1,10}$) \\
\hline
$1$ & $2$ & $3$ & 0.372 & \hs 0.957 & \hs ($B_{2 3}$) \\
\hline
$2$ & $2$ & $3$ & 0.564 & \hs 0.842 & \hs ($B_{2 3}$) \\
\hline
$2$ & $2$ & $10$ & 0.324 & \hs 0.955 & \hs ($B_{2,10}$) \\
\hline
$2$ & $3$ & $3$ & 0.419 & \hs 0.907 & \hs ($B_{3 3}$) \\
\hline
$2$ & $3$ & $4$ & 0.366 & \hs 0.948 & \hs ($B_{3 4}$) \\
\hline
$2$ & $3$ & $10$ & 0.284 & \hs 0.971 & \hs ($B_{3,10}$) \\
\hline
$2$ & $10$ & $10$ & 0.254 & \hs 0.909 & \hs ($B_{10,10}$) \\
\hline
$4$ & $5$ & $6$ & 0.417 & \hs 0.878 & \hs ($F$) \\
\hline
$10$ & $10$ & $11$ & 1.238 & \hs 0.817 & \hs ($F$) \\
\hline
$10$ & $11$ & $11$ & 0.969 & \hs 0.825 & \hs ($F$) \\
\hline
$10$ & $11$ & $12$ & 0.804 & \hs 0.845 & \hs ($F$) \\
\hline
$10$ & $11$ & $20$ & 0.391 & \hs 0.906 & \hs ($F$) \\
\hline
$10$ & $20$ & $20$ & 0.214 & \hs 0.888 & \hs ($F$) \\
\hline
\end{tabular}
\hskip 1 cm
\begin{tabular}{c|c|c||c|cl}
$\ell$ & $m$ & $n$ & ~~$K^{+}_{\ell m n d}~~ $ & $K^{-}_{\ell m n d} /K^{+}_{\ell m n d}$ \\
\hline \hline
$0$ & $1$ & $1$ & 0.135 & \hs 0.842 & \hs ($B_{2 2}$) \\
\hline
$0$ & $1$ & $2$ & 0.0694 & \hs 0.918 & \hs ($F$) \\
\hline
$0$ & $1$ & $10$ & 0.0215 & \hs 0.988 & \hs ($F$) \\
\hline
$1$ & $1$ & $2$ & $1/2 \sqrt{2 \pi}$ (${}^{\ast}$)
 & \hs 1 & \hs ($S$) \\
\hline
$1$ & $1$ & $3$ & 0.101 & \hs 0.987 & \hs ($S$) \\
\hline
$1$ & $1$ & $10$ & 0.0296 & \hs 0.995 & \hs ($S$) \\
\hline
$1$ & $2$ & $3$ & 0.0581 & \hs 0.865 & \hs ($F$) \\
\hline
$2$ & $2$ & $3$ & 0.115 & \hs 0.916 & \hs ($B_{2 3}$) \\
\hline
$2$ & $2$ & $10$ & 0.0302 & \hs 0.981 & \hs ($B_{2,10}$) \\
\hline
$2$ & $3$ & $3$ & 0.0646 & \hs 0.901 & \hs ($B_{3 3}$) \\
\hline
$2$ & $3$ & $4$ & 0.0482 & \hs 0.916 & \hs ($B_{3 4}$) \\
\hline
$2$ & $3$ & $10$ & 0.0237 & \hs 0.909 & \hs ($B_{3,10}$) \\
\hline
$2$ & $10$ & $10$ & 0.0167 & \hs 0.754 & \hs ($F$) \\
\hline
$4$ & $5$ & $6$ & 0.0437 & \hs 0.870 & \hs ($F$) \\
\hline
$10$ & $10$ & $11$ & 0.0990 & \hs 0.798 & \hs ($F$) \\
\hline
$10$ & $11$ & $11$ & 0.0734 & \hs 0.817 & \hs ($F$) \\
\hline
$10$ & $11$ & $12$ & 0.0583 & \hs 0.833 & \hs ($F$) \\
\hline
$10$ & $11$ & $20$ & 0.0223 & \hs 0.905 & \hs ($F$) \\
\hline
$10$ & $20$ & $20$ & 0.00978 & \hs 0.974 & \hs ($F$) \\
\hline
\end{tabular}
$$
}}
\vskip 1cm \noindent
\rule{6cm}{0.05mm}
\vskip 0.1cm \noindent
${~}^\ast$ {\footnotesize{Note that ${1 \over 2 \sqrt{2 \pi}}$ = 0.1994...~.
The equality $K^{-}_{1 1 2 3}/K^{+}_{1 1 2 3} = 1$ indicates
that ${1 \over 2 \sqrt{2 \pi}}$ is the sharp constant $K_{1 1 2 3}$.}}
\vfill \eject \noindent
\oddsidemargin=0.5truecm
\section{Asymptotics for the upper and lower bounds on $\boma{K_{\ell m n d}}$\,.}
\label{secasimp}
As reviewed in the Introduction, in our previous work on the constant
$K_{\ell \ell \ell d} \equiv K_{\ell d}$ we have analyzed the
$\ell \vain + \infty$ asymptotics of some upper and lower bounds for this constant,
the conclusion being \rref{eqintro}. \parn
Now, we are in condition to analyze more general limit cases;
here we discuss the behavior of $K_{\ell m n d}$ when
\beq m = b \,\ell, ~n = c \, \ell~ \quad  (1 \leqs b \leqs c), \qquad \ell \vain + \infty~. \feq
We note that conditions \rref{c1} on $\ell$, $m=b \, \ell$, $n= c \, \ell$ and $d$ are fulfilled if
\beq 1 \leqs b \leqs c~, \qquad \ell > {d \over 2(b + c - 1)}~. \label{bcell} \feq
Let us first analyze the asymptotics of an upper bound for $K_{\ell m n d}$.
Our starting point is the inequality
\beq K_{\ell m n d} \leqs K^{\SS}_{\ell m n d} := \sqrt{\sup_{u \in [0,+\infty)} \SS_{\ell m n d}(u)}~,
\label{eqks} \feq
with $\SS_{\ell m n d}$ as in Eq. \rref{ff}, to be used with $m= b \ell$
and $n = c \ell$. We note that Eqs. \rref{ff} \rref{fmnd} give
\beq \SS_{\ell, b \ell, c \ell, d}(u) =
{\Gamma((b + c) \ell - d/2) \over (4 \pi)^{d/2} \Gamma((b + c) \ell)}\,
\Sigma_{b c d \ell}(u)~, \label{eqss} \feq
\beq \Sigma_{b c d \ell}: [0,+\infty) \vain (0,+\infty)~, \label{eqsig} \feq
$$ u \mapsto \Sigma_{b c d \ell}(u) := (1 + 4 u)^\ell~
{~}_{3} F_{2}((b + c) \ell - {d \over 2}, b \ell, c \ell;
{(b + c) \ell \over 2}, {(b + c) \ell + 1 \over 2}; -u)~. $$
Our subsequent analysis rests on the condition
introduced hereafter.
\begin{prop}
\textbf{Definition.} Let $1 \leqs b\leqs c$, and $d \in \naturali_0$. We say that
\textsl{condition $\Co_{b c d}$ holds} if
\beq \sup_{u \in [0,+\infty)} \Sigma_{b c d \ell}(u) = 1 + O(1/\ell) \qquad \mbox{for $\ell
\vain + \infty$}~. \label{cobcd} \feq
\end{prop}
\begin{rema}
\textbf{Remarks.} (i) In any case, $\Sigma_{b c d \ell}(0) = 1$. So, the above
condition means that $\sup_{u} \Sigma_{b c d \ell}$ is close to the value
of the function at $u=0$. \parn
(ii) Condition $\Co_{1 1 d}$ \textsl{does not} hold for any $d \in \naturali_0$. In fact,
with the present notations, Proposition 2.2 of \cite{uno} gives
$\sup_{u \in [0,+\infty)} \Sigma_{1 1 d \ell}(u) = \Sigma_{1 1 d \ell}(1/2)[1 + O(1/\ell)]=
3^{d/2 + 1/2} 2^{-d/2} (4/3)^\ell [1 + O(1/\ell)]$ for $\ell \vain + \infty$. \parn
On the other hand, this negative result is not important for our purposes:
in fact the case $b=c=1$, i.e., $\ell = m = n$, is just the one analyzed by
different means in \cite{uno}, and summarized here
via Eq. \rref{eqintro}. \parn
Hereafter we consider a case where
$\Co_{b c d}$ can be proved, and another one where it can be reasonably conjectured.
\end{rema}
\begin{prop}
\label{co22}
\textbf{Proposition.} Condition $\Co_{2 2 d}$ holds for each $d \in \naturali_0$.
\end{prop}
\textbf{Proof.} See Section \ref{secbcd}. \fine
\begin{rema}
\label{remsigz}
\textbf{Remark.} The above result is sufficient for our purposes, but there is 
evidence for a slightly stronger statement:
$\sup_{u \geqs 0} \Sigma_{2 2 d \ell}$ is attained at
a point $u = u_{2 2 d \ell} \neq 0$ that, for $\ell \vain + \infty$, converges
to zero in such a way to fulfill condition \rref{cobcd}. We return to 
this point in the forthcomig Remark \ref{remsig}.
\fine
\end{rema}
Let us pass from the case $b=c=2$ to $b=2, c =3$; for the latter we have found
numerical evidence (but no analytic proof) for the following conjecture.
\begin{prop}
\label{conj}
\textbf{Conjecture.}
For each $d \in \naturali_0$ there is a real number $\ell_d > {d/8}$ such that, for all
$\ell \geqs \ell_d$, the function $\Sigma_{2 3 d \ell}$ is strictly decreasing on $[0,+\infty)$. So
\beq \sup_{u \in [0,+\infty)} \Sigma_{2 3 d \ell}(u) = \Sigma_{2 3 d \ell}(0) = 1  \qquad \mbox{for each
$\ell \geqs \ell_d$}~ \label{conjes} \feq
(which implies condition $\Co_{2 3 d}$, in a strong version with no term $O(1/\ell)$ in Eq.
\rref{cobcd}).
\end{prop}
In the above, the condition $\ell_d > d/8$ reflects the inequality on $\ell$ in Eq.
\rref{bcell}, for $b=2$ and $c=3$.
Conjecture \ref{conj} is probably related to
some inequalities for the ${}_{q+1} F_{q}$ functions, conjectured in \cite{KS}.
\begin{prop}
\label{ksas}
\textbf{Proposition.}
Suppose condition $\Co_{b c d}$ to hold for some fixed $b,c,d$ ($1 \leqs b \leqs c$, $d \in \naturali_0$).
Then, the upper bound $K^{\SS}_{\ell, b \ell, c \ell, d} $ on $K_{\ell, b \ell, c \ell, d}$ has the asymptotics
\beq K^{\SS}_{\ell, b \ell, c \ell, d} = {1 + O(1/\ell) \over
[4 (b + c) \pi \ell]^{d/4}} \quad \mbox{for $\ell \vain + \infty$}~. \label{kksas} \feq
\end{prop}
\textbf{Proof.} Let $\ell \vain + \infty$. Eqs. (\ref{eqks}-\ref{cobcd}) give
\beq K^{\SS}_{\ell, b \ell, c \ell, d} =
\sqrt{{\Gamma((b + c) \ell - d/2) \over \Gamma((b + c) \ell)}} \, {1 + O(1/\ell) \over (4 \pi)^{d/4}}~. \feq
Now, the thesis follows using the relation
\beq {\Gamma((b + c) \ell - d/2) \over \Gamma((b + c) \ell)} = {1 + O(1/\ell) \over {[(b+c) \ell}]^{d/2}}~, \feq
which is a consequence of Eq. \rref{rapgamma}. \fine
Let us pass to the asymptotics for a suitable lower bound on $K_{\ell, b \ell, c \ell, d}$.
We recall that, for any $\ell, m, n$, we have the Fourier lower bound \rref{theabovee};
let us use this with $p=q=0$. So, for all $\si, \tau \in (0,+\infty)$,
\beq K_{\ell m n d} \geqs \Kap^{F}_{\ell m n d}(\si, \tau) :=
{ \| f_{\si + \tau, d} \|_\ell \over \| f_{\si  d} \|_m \| f_{\tau d}\|_n}~;
\label{theabovee0} \feq
here $f_{\si d} := f_{p=0, \si, d}$, i.e.,
\beq f_{\si d} : \reali^d \vain \reali~, \qquad x \mapsto
f_{\si d}(x):= e^{-\si |x|^2/2} \qquad  (\,\sigma \in (0, +\infty)\,)~. \label{fsigd} \feq
Our main result in this framework is the following.
\begin{prop}
\label{kfas}
\textbf{Proposition.} Let $1 \leqs b \leqs c$, $d \in \naturali_0$, and
\beq \Delta_{b c} := \{ (\xi, \eta) \in (0,1/b) \times (0,1/c)~|~ \xi + \eta < 1~\}~. \feq
Then, for fixed $(\xi, \eta) \in \Delta_{b c}$ and $\ell \vain + \infty$,
\beq \Kap^{F}_{\ell,  b \ell, c \ell, d}\big({\xi \over \ell},
{\eta \over \ell}\big) = {1 + O(1/\ell) \over [D_{b c}(\xi,\eta) \pi \ell]^{d/4}}~,
\quad D_{b c}(\xi,\eta) := {(1 - \xi - \eta) (\xi + \eta) \over \xi \eta (1 - b \xi) (1 - c \eta)}~.
\label{relk} \feq
\end{prop}
\textbf{Proof.} See Section \ref{secbcd}. \fine
For given $b, c$ one uses Eq. \rref{relk} choosing $(\xi, \eta) \in \Delta_{b c}$
so as to minimize $D_{b c}$ (or to go
as close as possible to the minimum point of this function); this choice gives the best lower bound
of the type \rref{relk}, in the limit $\ell \vain + \infty$. \parn
Let us write down two Corollaries of Propositions \ref{ksas} and \ref{kfas}, for
the cases $b=c=2$ and $b=2, c=3$, respectively.
\begin{prop}
\label{cor22}
\textbf{Corollary}. For any $d \in \naturali_0$, the following holds. \parn
(i) The upper bound $K^{\SS}_{\ell, 2 \ell, 2 \ell, d}$ is such that
\beq K^{\SS}_{\ell, 2 \ell, 2 \ell, d} = {1 + O(1/\ell) \over
(16 \pi \ell)^{d/4}} \qquad \mbox{for $\ell \vain + \infty$}~. \feq
(ii) The function $D_{2 2} : \Delta_{2 2} \vain (0,+\infty)$ from Proposition \ref{kfas}
attains its minimum at $\xi=\eta=1/4$. It is $D_{2 2}(1/4,1/4) = 16$; so,
the corresponding lower bound $\Kap^{F}_{\ell,  2 \ell, 2 \ell, d}\big((1/(4 \ell),
1/(4 \ell)\big)$ is such that
\beq \Kap^{F}_{\ell,  2 \ell, 2 \ell, d}\big({1 \over 4 \ell},
{1 \over 4 \ell}\big) = {1 + O(1/\ell) \over (16 \pi  \ell)^{d/4}}
\qquad \mbox{for $\ell \vain + \infty$}~. \feq
(iii) As a consequence of (i) (ii), one has
\beq  K_{\ell, 2 \ell, 2 \ell, d} = {1 + O(1/\ell) \over
(16 \pi \ell)^{d/4}} \qquad \mbox{for $\ell \vain + \infty$}~. \feq
\end{prop}
\textbf{Proof.} (i) Use Proposition \ref{ksas} with $b=c=2$ (recalling that condition $\Co_{2 2}$ holds, by
Proposition \ref{co22}). \parn
(ii) Elementary. \parn
(iii) The thesis follows from $\Kap^{F}_{\ell,  2 \ell, 2 \ell, d}\big({1 \over 4 \ell},
{1 \over 4 \ell}\big) \leqs K_{\ell, 2 \ell, 2 \ell, d} \leqs K^{\SS}_{\ell, 2 \ell, 2 \ell, d}$. \fine
\begin{prop}
\label{cor23}
\textbf{Corollary}. For $d \in \naturali_0$, we have the following. \parn
(i) If $\Co_{2 3 d}$ holds, the upper bound $K^{\SS}_{\ell, 2 \ell, 3 \ell, d}$ is such that
\beq K^{\SS}_{\ell, 2 \ell, 3 \ell, d} = {1 + O(1/\ell) \over
(20 \pi \ell)^{d/4}} \qquad \mbox{for $\ell \vain + \infty$}~. \feq
(ii) Consider the function $D_{2 3} : \Delta_{2 3} \vain (0,+\infty)$ from Proposition \ref{kfas},
and evaluate it at $(\xi, \eta) := (1/5, 1/7)$
(which is close to its minimum point). It is $D_{2 3}(1/5,1/7) = 23$; so
the corresponding lower bound $\Kap^{F}_{\ell,  2 \ell, 3 \ell, d}\big(1/(5 \ell),
1 / (7 \ell)\big)$ is such that
\beq \Kap^{F}_{\ell,  2 \ell, 3 \ell, d}\big({1 \over 5 \ell},
{1 \over 7 \ell}\big) = {1 + O(1/\ell) \over (23 \pi  \ell)^{d/4}}
\qquad \mbox{for $\ell \vain + \infty$}~. \feq
(iii) Summing up, (i) (ii) give
\beq  {1 + O(1/\ell) \over (23 \pi  \ell)^{d/4}} \leqs
K_{\ell, 2 \ell, 2 \ell, d} \mpe {1 + O(1/\ell) \over
(20 \pi \ell)^{d/4}} \qquad \mbox{for $\ell \vain + \infty$}~, \feq
\vskip -0.9cm \noindent
where $\mpe$ means that the indicated relation is true
if condition $\Co_{2 3 d}$ holds.
\end{prop}
\textbf{Proof.} (i) Use Proposition \ref{ksas} with $b=2, c=3$.
(ii) Elementary. (iii) Obvious. \fine
\section{Proof of Proposition \ref{pupper}.}
\label{upper}
Here and in the rest of the paper, we work in a fixed
dimension $d \in \naturali_0$.
The proof of
the cited proposition is preceded by some lemmas.
The method is similar to the one of \cite{uno}, but technically more difficult; again,
the basic idea is to work with the Fourier transform $\FF$, that sends the pointwise product
of functions into the convolution. \parn
Let us write $F \ast G$ for the convolution of two
complex functions $F, G$ on $\reali^d$,
given by
\beq (F \ast G)(k) := \int_{\reali^d}
d h~F(k - h) G(h)~.  \feq
We have
\beq \FF (f g) = {1 \over (2 \pi)^{d/2}}~\FF f \ast \FF g \label{send} \feq
for all sufficiently regular functions
$f$ and $g$ on $\reali^d$ (and, in particular, for functions
to which we will apply \rref{send} in the rest of the section).
\parn
Let us recall the definition \rref{efnd}
$\F_{t d}(k) := 1/(1 + | k |^2)^t$ for all $t \in \reali$ and $k \in \reali^d$,
to which we will refer systematically in the sequel.
The forthcoming Lemmas consider pairs $m, n$ or triples $\ell, m, n$
of real numbers.
\begin{prop}
\label{pupperz}
\textbf{Lemma.} Let $m + n > d/2$.
Then, the integral defining the convolution $(\F_{m d} \ast \F_{n d})(k)$
is convergent, for all $k \in \reali^d$.
\end{prop}
\textbf{Proof.} For an integral $\int_{\reali^d} F(h) d h$ to be
convergent, it suffices that $F$ be continuous and
that, for $h \vain \infty$, $F(h) = O(1/|h|^\eta)$ with $\eta > d$.
For any $k \in \reali^d$, the convolution integral
\beq (\F_{m d} \ast \F_{n d})(k) = \int_{\reali^d} {d h \over (1 + |k -h|^2)^m (1 + |h|^2)^n}~
\label{convint} \feq
fulfills these conditions with $\eta = 2(m + n)$. \fine
\begin{prop}
\label{pupper1}
\textbf{Lemma.} Let $\ell, m, n$ fulfill \rref{c1}. Then
\beq K_{\ell m n d} \leqs \sqrt{  \sup_{k \in \reali^d} \ES_{\ell m n d}(k) }~, \feq
\beq \ES_{\ell m n d}(k) := {(1 + | k |^2)^\ell \over (2 \pi)^d}
\left(\F_{m d} \ast \F_{n d} \right)(k)~. \label{esnd} \feq
\end{prop}
\textbf{Proof.} Consider any two functions $f \in \Hm$, $g \in \Hn$. Then
\beq \| f g \|_\ell^2 = \int_{\reali^d} d k (1 + k^2)^\ell | \FF(f g)(k) |^2 = {1 \over (2 \pi)^{d}}
\int_{\reali^d} d k (1 + k^2)^\ell | (\FF f \ast  \FF g)(k) |^2~. \label{laprec} \feq
Explicitating the convolution we find
\beq  (\FF f \ast  \FF g)(k)  =  \int_{\reali^d} \!\!\!dh\, \FF f(k - h) \FF g (h)  \feq
$$ = \int_{\reali^d} \!\!\! dh {1 \over \sqrt{1 + |k - h|^2}^{\,m} \sqrt{1 + | h |^2}^{\,n}}~
( \sqrt{1 + |k - h|^2}^{\,m}~\FF f(k - h) \sqrt{1 + |h|^2}^{\,n} \FF g (h)), $$
and H\"older's inequality $| \int d h~ U(h) V(h) |^2 \leqs \int d h | U(h) |^2 \,
\int d h~| V(h) |^2 $ gives
\parn
\vbox{
\beq | (\FF f \ast  \FF g)(k) |^2 \leqs C(k) P(k)~, \label{dains} \feq
$$ C(k) := \int_{\reali^d} {d h \over (1 + | k - h |^2)^m (1 + | h |^2)^n } =
\left(\F_{m d} \ast \F_{n d}\right)(k)~, $$
$$ P(k) := \int_{\reali^d} d h (1 + |k - h|^2)^m | \FF f (k - h) |^2
(1 + | h |^2)^n | \FF g(h) |^2~. $$ }
Inserting \rref{dains} into Eq. \rref{laprec} we get
\beq \| f g \|_\ell^2 \leqs {1 \over (2 \pi)^{d}} \int_{\reali^d} d k (1 + | k |^2)^\ell C(k) P(k) \feq
$$ \leqs \Big(\sup_{k \in \reali^d} {(1 + | k |^2)^\ell \over (2 \pi)^d} C(k) \Big)~\int_{\reali^d} dk \, P(k)~
= \Big(\sup_{k \in \reali^d} \ES_{\ell m n d}(k) \Big) \int_{\reali^d} dk \, P(k)~. $$
But
$$\int_{\reali^d} d k \, P(k) =  \int_{\reali^d} d k (1 + |k |^2)^m | \FF f(k) |^2 ~
\int_{\reali^d} d h (1 + |h |^2)^n | \FF g(h) |^2 = \| f \|_m^2 \, \| g \|_n^2~, $$
so we are led to the thesis. \fine
\begin{prop}
\label{pupper2a}
\textbf{Lemma.} Let $m, n \geqs 0$, $m + n > d/2$. Then,
for all $k \in \reali^d$,
\beq \left(\F_{m d} \ast \F_{n d} \right)(k) =
\pi^{d/2} {\Gamma(m + n - d/2) \over \Gamma(m + n)}
F_{m n d}\big({| k |^2 \over 4}\big)~, \label{cleara} \feq
where $F_{m n d}$ is the hypergeometric function (of the
${}_3 F_{2}$ type) in Eq. \rref{fmnd} of Proposition \ref{pupper}.
\end{prop}
\textbf{Proof.}
Both sides of \rref{cleara} are symmetric in $m, n$, so we can
restrict the attention to the case $m \leqs n$ and write our
basic assumptions as
\beq 0 \leqs m \leqs n,~~~m + n > {d \over 2}~. \label{c1mn} \feq
Conditions \rref{c1mn} on $m, n$ are equivalent to
\beq {d \over 4} < n,~ m \in M_{n d}~,
\qquad M_{n d} := [0, n] \cap (d/2 -n, + \infty)~. \label{cc1mn} \feq
Let us fix $k \in \reali^d$. We claim that it is sufficient
to prove the thesis \rref{cleara} under even more restrictive conditions than \rref{cc1mn},
namely, for
\beq {d \over 4} < m \leqs n~. \label{cc2mn} \feq
In fact, for fixed ($k \in \reali^d$ and) $n > d/4$: \parn
(i) both sides of Eq.\,\rref{cleara}, viewed as functions of $m$, are analytic in an
open neighborhood on $M_{n d}$, namely, the interval $(d/2 - n, + \infty)$. This is
made evident by the expression \rref{convint} for the
convolution integral $\left(\F_{m d} \ast \F_{n d} \right)(k)$
and by the considerations about ${}_{q+1} F_{q}$ following Eq. \rref{intef},
here applied to $F_{m n d}(|k|^2/4) = {}_3 F_2 (m + n - d/2, m, n; (m + n)/2, (m + n + 1)/2,-|k|^2/4)$
({\footnote{The analyticity result for ${}_{q+1} F_{q}$ stated
after the integral representation \rref{intef} ensures the
following in the present case, for fixed $n > d/4$ and $u \in [0,+\infty)$:
the function $m \mapsto {}_3 F_2 (m + n - d/2, m, n; m/2 + n/2, m/2 + n/2 + 1/2,-u)$
is analytic where $m$ fulfills the condition $m + n - d/2 > 0$, i.e., for
$m \in (d/2 - n, + \infty)$. }}). \parn
(ii) By the principle of analytic continuation, if the
two sides of \rref{cleara} are equal for $m \in (d/4, n]$, they
are equal as well for $m$ in ${M}_{n d}$. \vskip 0.1cm \noindent
The rest of the proof is devoted to establishing \rref{cleara} for
$m, n$ as in \rref{cc2mn}. \parn
Under these conditions we can represent $\F_{t d}$ as the Fourier transform
of the function $\f_{t d}$ (Eqs. \rref{fnd} \rref{gemac}), both
for $t=n$ and for $t =m$. From here and
\rref{send},
\beq \left(\F_{m d} \ast \F_{n d} \right)(k)
= (2 \pi)^{d/2} \FF\left(\f_{m d} \f_{n d} \right)(k)~. \label{es} \feq
The product $\f_{m d} \f_{n d}$ is a radially symmetric function, whose explicit expression in terms of
Macdonald functions follows from \rref{gemac}. So,
$\FF\left(\f_{m d} \f_{n d} \right)$ can be computed using the formula
\rref{eboc} for radially symmetric Fourier transforms, and the conclusion is
\beq \left(\F_{m d} \ast \F_{n d} \right)(k) \label{brf} \feq
$$ = {(2 \pi)^{d/2} \over {2^{m + n - 2} \Gamma(n) \Gamma(m) | k |^{d/2 - 1}}}~
\int_{0}^{+\infty} \!\!\! dr \, r^{m + n - d/2} J_{d/2 - 1}(| k | r)~K_{m - d/2}(r) K_{n - d/2}(r)~; $$
the above integral is computed via \rref{result}, and in this way one gets the thesis
\rref{cleara}. \parn
(Final remark: some of our last
manipulations seem to exclude the point $k=0$, see e.g. the denominator in Eq. \rref{brf}; however,
Eq. \rref{cleara} holds here as well, by continuity). \fine
\vskip 0.2cm \noindent
\begin{prop}
\label{pupper2}
\textbf{Lemma.} Let $\ell, m, n$ fulfill \rref{c1}. Then, for all $k \in \reali^d$,
\beq \ES_{\ell m n d}(k) = \SS_{\ell m n d}\,({| k |^2 / 4})~, \label{clear} \feq
where $\SS_{\ell m n d}$ is the function in Eq. \rref{ff} of Proposition \ref{pupper}.
\end{prop}
\textbf{Proof.} This follows immediately from the definition
\rref{esnd}
$\ES_{\ell m n d}(k)$ $:= \dd{(1 + | k |^2)^\ell \over (2 \pi)^d}$
$\left(\F_{m d} \ast \F_{n d} \right)(k)$, from Eq. \rref{cleara}
of the previous Lemma and from the definition \rref{ff} of
$\SS_{\ell m n d}$.
\fine
We are finally ready to derive the main result of the section, i.e.,
to prove Proposition \ref{pupper}.
\vskip 0.2cm \noindent
\textbf{Proof of Proposition \ref{pupper}, item (i).} Again,
$\ell, m, n$ are assumed to fulfill \rref{c1}. Lemmas \ref{pupper1} and \ref{pupper2} give
immediately the bound \rref{kpnd} for $K_{\ell m n d}$, with
$\SS_{\ell m n d}$ as in Eq. \rref{ff}; in the sequel we frequently
mention the hypergeometric function $F_{m n d}$ appearing in Eqs. \rref{ff} \rref{fmnd},
recalling again that this is of the ${~}_3 F_{2}$ type. \parn
In the special case $m=n$,
the expression \rref{fmmd} of $F_{m n d}$ as a ${~}_2 F_{1}$ function
follows immediately from \rref{ofco}. Eq. \rref{ftriv} for the ''trivial'' case
$m=0$ arises noting that
$F_{0 n d}(u) = {}_3 F_2(n - d/2, 0, n; n/2, n/2 + 1/2, -u) = 1$
by \rref{limitriv}. \parn
Let us prove the properties of
$\SS_{\ell m n d}$ mentioned in item (i), for arbitrary $\ell, m, n, d$. \parn
First of all, the statement $\SS_{\ell m n d}(u) \in (0,+\infty)$ for all $u \in [0,+\infty)$
follows immediately from the relation \rref{clear} between this function and
$\ES_{\ell m n d}$, which is positive due to the definition \rref{esnd}. Any hypergeometric
function ${}_p F_q$ takes the value $1$ at the origin; so, $\SS_{\ell m n d}(0)$ has the expression
\rref{duesei}. To conclude, we must prove the asymptotics \rref{dueset} for $\SS_{\ell m n d}(u)$
as $u \vain + \infty$; this will give the result \rref{dueot} for
$\SS_{\ell m n d}(+\infty)$, also implying the boundedness of $\SS_{\ell m n d}$ on $[0,+\infty)$.
\parn
To derive \rref{dueset}, we first consider the case $m < n$ and apply to $F_{m n d}(u)$
the general asymptotics \rref{asimp}
(with $\alpha=m$, $\beta=n$,
$\gamma=m + n - d/2$); with the obvious relation $(1 + 4 u)^\ell \sim (4 u)^\ell$, this gives
\beq \SS_{\ell m n d}(u) \sim {4^\ell \over (4 \pi)^{d/2}}
{\Gamma(n - {d \over 2}) \over \Gamma(n)} {\Gamma_{m n} \over u^{m - \ell}}
\qquad \mbox{for $u \vain + \infty$}~, \label{ss2} \feq
$$ \Gamma_{m n} :=  {\Gamma({m \over 2} + {n \over 2})
 \Gamma({m \over 2} + {n \over 2} + {1 \over 2})
\over \Gamma(m + n)} {\Gamma(n-m) \over \Gamma({n \over 2} - {m \over 2})
\Gamma({n \over 2} - {m \over 2} + {1 \over 2})}~. $$
On the other hand, expressing $\Gamma(n \pm m)$ via the duplication formula \rref{dupl}
we see that
\beq \Gamma_{m n} = {1 \over 4^m} \quad \mbox{for all $n$} \label{ss1}~; \feq
Eqs. \rref{ss1} and \rref{ss2} give the thesis \rref{dueset}, with the previous assumption $m < n$.
To conclude, we must derive
\rref{dueset} in the special case $m=n$, where $F_{m n d}$ collapses into a ${~}_2 F_{1}$ function
due to \rref{fmmd}; this case is worked out similarly to the previous one, using
the asymptotics \rref{asi21} (and again, the duplication formula for $\Gamma$).
\fine
\textbf{Proof of Proposition \ref{pupper}, item (ii).} Our aim
is to derive the series expansions for $F_{m n d}$ in the cited item of the proposition,
and to show that they are just finite sums with the special assumptions on $m,n,d$ indicated therein. \parn
First of all we note that, for $u \in [0,+\infty)$,
\parn
\vbox{
\beq F_{m n d}(u) \label{espp} \feq
$$ = \sum_{i=0}^{+\infty} { \left(m + n - {d \over 2} \right)_{i} (m)_i
\left( {m - n + 1 \over 2} \right)_{i} \over
\left({m + n  \over 2} \right)_{i} \left({m + n + 1 \over 2} \right)_{i}  }~{u^i \over i!}
{~}_{2} F_{1}(m + n - {d \over 2} + i, m + i; {m + n  \over 2} + i; -u) $$
$$ = \sum_{i=0}^{+\infty} { \left(m + n - {d \over 2} \right)_{i} (m)_i
\left( {m - n \over 2} \right)_{i} \over \left({m + n + 1 \over 2}
\right)_{i} \left({m + n  \over 2} \right)_{i} }~{u^i \over i!}
{~}_{2} F_{1}(m + n - {d \over 2} + i, m + i; {m + n + 1 \over 2} + i; -u)~. $$}
In the above, the first equality follows directly from the definition
\rref{fmnd} and from the expansion \rref{dadare}; the
second equality follows writing $F_{m n d}(u) =
{~}_{3} F_{2}(m + n - {d \over 2}, m, n; {m + n + 1 \over 2}, {m + n \over 2} ; -u)$,
and then using again Eq. \rref{dadare}.
On the other hand,
$$ {~}_{2} F_{1}(m + n - {d \over 2} + i, m + i; {m + n  \over 2} + i; -u)
= { {~}_{2} F_{1}(\dd{d - m - n \over 2}, \dd{n - m \over 2}; \dd{m + n \over 2} + i; -u)
\over (1 + u)^{{3 m + n - d \over 2} + i}} $$
\beq = {1 \over (1 + u)^{{3 m + n - d \over 2} + i}}
\sum_{j=0}^{+\infty} {\left( {d - m - n \over 2} \right)_j
\left( {n - m \over 2} \right)_j \over \left( {m + n \over 2} + i \right)_j } {(-u)^j \over j!} ~;
\label{211} \feq
the first equality above follows from the Kummer transformation \rref{sukum},
the second one reflects the standard power series expansion \rref{partpq} for ${~}_{2} F_{1}$.
The latter expansion holds if $u \in [0,1)$, or $u \in [0,+\infty)$ and the series
over $j$ is a finite sum; these are just the conditions in the Proposition under proof. \parn
Inserting the expansion \rref{211} into the first equality \rref{espp}, one gets
\rref{esp1}. \parn
For similar reasons, we can write \parn
\vbox{
$$ {~}_{2} F_{1}(m + n - {d \over 2} + i, m + i; {m + n + 1 \over 2} + i; -u)
= { {~}_{2} F_{1}({d + 1 - m - n \over 2}, {n + 1 - m \over 2}; {m + n + 1 \over 2} + i; -u)
\over (1 + u)^{{3 m + n - d - 1\over 2} + i}}
$$
\beq = {1 \over (1 + u)^{{3 m + n - d - 1\over 2} + i}}
\sum_{j=0}^{+\infty} {\left( {d + 1 - m - n \over 2} \right)_j
\left( {n + 1 - m \over 2} \right)_j \over \left( {m + n + 1 \over 2} + i \right)_j } {(-u)^j \over j!}  \feq}
(again when $u \in [0,1)$, or $u \in [0,+\infty)$ and the series over $j$ is a finite sum).
Inserting this result into the second equality \rref{espp}, one gets
\rref{esp2}. \parn
We finally come to statements (\ref{spe1}-\ref{spe2}), giving conditions
for the series over $j, i$ in \rref{esp1} or \rref{esp2}
to become finite sums; as an example, we account for the first of such statements.
\parn
The series over $j$ in \rref{esp1} contains the Pochhammer symbol
$\left({d - m - n \over 2}\right)_j$; on the other hand, the
assumption in the first line of \rref{spe1} is equivalent to
\beq {d - m - n \over 2} = - h, \quad h \in \naturali~. \feq
From $h \in \naturali$ we infer $(-h)_j = 0$ for $j > h$, so
\beq  \sum_{j=0}^{+\infty} \vain \sum_{j=0}^{h} = \sum_{j=0}^{m + n - d \over 2}
\quad \mbox{in \rref{esp1}}~. \feq
The other statements in (\ref{spe1}-\ref{spe2}) are proved analyzing:
the term $\left({m - n + 1 \over 2}\right)_i$ in \rref{esp1};
the term $\left({d + 1 - m - n  \over 2}\right)_j$ in \rref{esp2};
the term $\left({m - n \over 2}\right)_i$ in \rref{esp2}.
\fine
\section{Proof of Proposition \ref{pbessel}.}
\label{pbes}
Hereafter we prove items (i) (ii) of the cited proposition (after this,
item (iii) will be obvious).
\vskip 0.2cm \noindent
(i) We must show that $\f_{\nu t d}$ belongs to $H^n(\reali^d, \complessi)$, and justify
the expression \rref{fnutd} for its $H^n$ norm. The relation $\f_{\nu t d} \in H^n(\reali^d, \complessi)$
follows from the finiteness of the integrals appearing below; the norm of this function is
given by
\beq \| \f_{\nu t d} \|^2_n = \int_{\reali^d} d k~(1 + | k |^2)^n
| \FF \f_{\nu t d}(k) |^2 = {1 \over \nu^{2 d}}~
\int_{\reali^d} d k~{(1 + | k |^2)^n \over (1 + | k |^2/\nu^2)^{2 t}}
\label{thelast} \feq
$$ = {2 \pi^{d/2} \over \Gamma(d/2) \nu^{2 d}}~
\int_{0}^{+\infty} \!\!\!\! d \rho \, \rho^{d-1}~{(1 + \rho^2)^n \over (1 + \rho^2/\nu^2)^{2 t}} =
{\pi^{d/2} \over \Gamma(d/2) \nu^{d}}~
\int_{0}^{+\infty} \!\!\!\! d u \, u^{d/2-1}~{(1 + \nu^2 u)^n \over (1 + u)^{2 t}}~. $$
In the last two passages we have used Eq. \rref{rag} for the integral of a radially symmetric function,
depending only on
$\rho := | k |$, and then we have changed the
variable $\rho$ to $u = \rho^2/\nu^2$. \parn
Let us fix the attention to the integral over $u$
(clearly convergent, due to the assumption $t > n/2 + d/4$ in the statement
under proof);
this integral is computed via the identity \rref{gf}, and one gets the
thesis \rref{fnutd}. \parn
(ii) In the proof of Lemma \ref{pupper2}, we have derived Eq. \rref{brf} for a Fourier transform of the type
$\FF (\f_{m d} \f_{n d})$. With similar manipulations, in this case we get
\beq \FF\left(\f_{\mu s d} \f_{\nu t d}\right)(k) \label{brft} \feq
$$ = {\mu^{s - d/2} \nu^{\,t -d/2} \over {2^{s + t - 2} \Gamma(s) \Gamma(t) | k |^{d/2 - 1}}}~
\int_{0}^{+\infty} \!\!\! dr~ r^{s + t - d/2} J_{d/2 - 1}(| k | r)~K_{s - d/2}(\mu r) K_{t - d/2}(\nu r)~, $$
and a coordinate change $r \vain r/|k|$ gives
\beq \FF\left(\f_{\mu s d} \f_{\nu t d}\right)(k) = G_{s t d}(\mu, \nu; {|k|^2 / 4})~, \feq
with $G_{s t d}$ as in \rref{gstd}. This implies
\parn
\vbox{
\beq \| \f_{\mu s d} \f_{\nu t d} \|^2_{\ell} = \int_{\reali^d} d k \,(1 + |k|^2)^\ell\,
|\FF\left(\f_{\mu s d} \f_{\nu t d}\right)(k)|^2 \label{into} \feq
$$= \int_{\reali^d} d k \, (1 + |k|^2)^\ell \, G_{s t d}^2(\mu, \nu; {|k|^2 / 4})~. $$}
On the other hand, for radial integrals we have $dk = 2 \, \pi^{d/2}
|k|^{d-1} \, d|k | \, /\Gamma(d/2)$, and putting $|k| = 2 \sqrt{u}$ we get
the expression \rref{fmufnu}  for $\| \f_{\mu s d} \f_{\nu t d} \|^2_{\ell}$. \parn
Finally, let us consider the case
$s - {d \over 2}, t - {d \over 2} \in \naturali + {1 \over 2}, ~~\ell \in \naturali$, and
show that Eqs. \rref{fmufnu} \rref{gstd} yield Eq. \rref{pooint}.
To this purpose, we first compute the function $G_{s t d}(\mu, \nu; u)$ in \rref{gstd};
in this case Eq. \rref{cita} for the Macdonald functions gives
\parn
\vbox{
\beq G_{s t d}( \mu, \nu; u) =
{\pi \over {2^{2 s + 2 t - 2} \Gamma(s) \Gamma(t)  }}   \label{gstdd} \feq
$$ \times \sum_{i=0}^{~s-{d \over 2} - {1 \over 2}} \sum_{j=0}^{~t- {d \over 2} - {1 \over 2}}
{(2 s - i - d - 1)! \, (2 t - j - d - 1)! \, \mu^i \nu^j \over i! \, j!
(s - i - {d \over 2} - {1 \over 2})! \, (t - j - {d \over 2} - {1 \over 2})!
\, u^{i/2 + j/2 + d/2}} $$
$$ \times \int_{0}^{+\infty}  \! \! \! \!d r~ r^{i + j + d/2} J_{d/2 - 1}(r)~e^{- {(\mu + \nu) r \over 2 \sqrt{u}}}~. $$}
On the other hand, for any $\sigma \in (0,+\infty)$,
\parn
\vbox{
\beq \int_{0}^{+\infty}  \! \! \! \!d r~ r^{i + j + d/2} J_{d/2 - 1}(r) e^{-r/\sigma} \label{su1} \feq
$$ = {(i+j + d-1)! \, \sigma^{i + j + d} \over 2^{d/2 - 1} \Gamma(d/2)}{~}_2 F_1 ({i + j + d \over 2}, {i + j + d + 1 \over 2};
{d \over 2}; - \sigma^2)  $$
$$ = {(i+j + d-1)! \, \sigma^{i + j + d} \over 2^{d/2 - 1}  \Gamma(d/2)(1 + \si^2)^{i + j + d/2 + 1/2}}
{~}_2 F_1 (-{i + j \over 2}, -{i + j + 1 \over 2}; {d \over 2}; - \sigma^2)~, $$}
where the first equality follows from \cite{Wat} (page 385, Eq. (2)), and the second one from the Kummer transformation
\rref{sukum}. Since $i, j$ are nonnegative integers, one of the numbers
${i + j \over 2}$ and ${i + j + 1 \over 2}$ is a nonnegative integer and equals $[ {i + j + 1 \over 2}]$; so,
\beq {~}_2 F_1 (-{i + j \over 2}, -{i + j + 1 \over 2}; {d \over 2}; - \sigma^2) =
\sum_{k=0}^{[ {i + j + 1 \over 2}]} {\left(-{i + j \over 2}\right)_k \left(-{i + j + 1 \over 2}\right)_k
\over \left({d \over 2}\right)_k} {(-1)^k \sigma^{2 k} \over k!}~. \label{su2} \feq
Now, setting $\si := 2 \sqrt{u}/(\mu + \nu)$ we substitute \rref{su2} into \rref{su1} and
then put the result into \rref{gstdd}; the conclusion is
\parn
\vbox{
\beq G_{s t d}( \mu, \nu; u) = {\pi \over \Gamma(d/2) \Gamma(s) \Gamma(t)}
\sum_{i=0}^{~s-{d \over 2} - {1 \over 2}} \sum_{j=0}^{~t- {d \over 2} - {1 \over 2}}
\sum_{k=0}^{[ {i + j + 1 \over 2}]}  \label{poi} \feq
$$ \times {(-1)^k (i + j + d-1)! (2 s - i - d - 1)! (2 t - j - d - 1)!
\left(-{i + j \over 2}\right)_k \left(-{i + j + 1 \over 2}\right)_k
\over 2^{2 s + 2 t - i - j - 2 k - d/2 - 3} \, i! \, j! \, k!
(s - i - {d \over 2} - {1 \over 2})! \, (t - j - {d \over 2} - {1 \over 2})!
\left({d \over 2}\right)_k}  $$
$$ \times {\mu^i \nu^j (\mu + \nu)^{i + j - 2 k + 1} u^{k} \over
\big( (\mu + \nu)^2 + 4 u \big)^{i + j + d/2 + 1/2}}~. $$}
The result \rref{poi} has the form
\beq G_{s t d}( \mu, \nu; u) = {\pi \over \Gamma(d/2) \Gamma(s) \Gamma(t)}
\sum_{(i j k) \in I_{s t d}}  G_{s t i j k d}
{\mu^i \nu^j (\mu + \nu)^{i + j - 2 k + 1} (4 u)^{k} \over
\big( (\mu + \nu)^2 + 4 u \big)^{i + j + d/2 + 1/2}}~, \feq
where $I_{s t d}$ and $G_{s t i j k d}$ are as in Eqs. \rref{istd} \rref{gstijkd}.
The next step is to insert this result into Eq.
\rref{fmufnu}  for $\| \f_{\mu s d} \f_{\nu t d} \|^2_{\ell}$; this contains the integral over $u$
of the expression \parn
\vbox{
\beq (1 + 4 u)^{\ell} G_{s t d}^2( \mu, \nu; u) = {\pi^2 \over \Gamma^2(d/2) \Gamma^2(s) \Gamma^2(t)}
\sum_{h=0}^{\ell} \left( \barray{c}  \ell \\ h \farray \right)
(4 u)^h   \feq
$$ \times \sum_{(i j k) \in I_{s t d}} \sum_{(i' j' k') \in I_{s t d}} G_{s t i j k d}
G_{s t i' j' k' d}  $$
$$ \times \, {\mu^{i+i'} \nu^{j+j'} (\mu + \nu)^{i + i' + j + j' - 2 k - 2 k' + 2}
(4 u)^{k + k'} \over
\big( (\mu + \nu)^2 + 4 u \big)^{i + i' + j + j' + d + 1}}~; $$}
we substitute this in \rref{fmufnu} and integrate over $u$, taking into account that
\beq \int_{0}^{+\infty} \! \! \! \! d u \, {(4 u)^{a} \over ( \xi + 4 u)^b} =
{\Gamma(a+1) \Gamma(b-a-1) \over 4 \, \Gamma(b)\,\xi^{b-a-1}}~. \feq
The conclusion is Eq. \rref{pooint} for $\| \f_{\mu s d} \f_{\nu t d} \|^2_{\ell}$. \fine
\section{Proof of Proposition \ref{propdir}.}
\label{secdir}
Throughout the section we make the assumptions of  Eq. \rref{lett}:
$$ 0 \leqs \ell \leqs n~, \qquad n > {d \over 2}~. $$
\begin{prop}
\textbf{Lemma.} One has
\beq K_{\ell \ell n d} \geqs {| g(0) | \over ~\| g \|_n} \label{tesilem} \feq
for each nonzero $g \in \Hn$. (Note that $g(0)$ makes sense, by the
well known imbedding $\Hn \subset C(\reali^d, \complessi)$.)
\end{prop}
\textbf{Proof.} Let us present the idea heuristically.
We fix a nonzero $g \in \Hn$, and write the inequality
\beq K_{\ell \ell n d} \geqs {\| f_{\ep} g \|_\ell \over \| f_{\ep} \|_{\ell} \| g \|_n} \label{fep} \feq
where $(f_{\ep})_{\ep > 0}$ is a family of approximants of the Dirac $\delta$ distribution on
$\reali^d$: $f_{\ep} \vain \delta$ as $\ep \vain 0^{+}$. Then, for $\ep \vain 0^{+}$,
$f_{\ep} g \sim  g(0) f_{\ep}$ and
\beq \| f g \|_\ell \sim | g(0) |~ \| f_{\ep} \|_\ell~; \feq
so, in this limit, the inequality \rref{fep} gives the thesis \rref{tesilem}.
For a rigorization of this argument, see the the proof of Lemma 7.1 in \cite{mp2}
(which contains a statement very similar to the present one). \fine
\textbf{Proof of Proposition \ref{propdir}.} From the previous Lemma,
\beq K_{\ell \ell  n d} \geqs \sup_{g \in \Hn \setminus \{0\}} {| g(0) | \over ~\| g \|_n}~; \feq
as shown in \cite{imb}, the above sup equals $S_{\infty n d}$ (and is attained at $g = g_{n d}$ as in Eqs.
\rref{fnd} \rref{gemac}). \fine
\section{Proof of Propositions \ref{co22} and \ref{kfas}.}
\label{secbcd}
Each one of the two proofs will be preceded by a lemma about the asymptotics of a Laplace integral;
we use this expression
to indicate an integral of the form
\beq L(\lambda) := \int_{0}^b dt \, \theta(t) \, e^{-\lambda \varphi(t)}
\qquad \big(b \in [0,+\infty), \lambda \in (\lambda_0, + \infty)~\big) \feq
where $\theta \in C((0,b),\reali)$, $\varphi \in C([0,b), \reali) \cap C^1((0,b), \reali)$
are such that $\int_{0}^b d t |\theta(t)| e^{-\lambda \varphi(t)}$ $< + \infty$ for all $\lambda$
as above, and $\varphi'(t) > 0$ for $t \in (0,b)$ (the prime meaning $d/ d t$).
The following implication is well known (see e.g. \cite{Olv}):
$$ {\theta(t) \over \varphi'(t)} = P (\varphi(t) - \varphi(0))^{\alpha -1}
[ 1 + O(\varphi(t)- \varphi(0))]~~\mbox{for $t \vain 0^{+}$} \quad \big(P
\in \reali, \alpha  \in (0,+\infty)~\big) $$
\beq \Longrightarrow L(\lambda) = P e^{-\lambda \varphi(0)}
{\Gamma(\alpha) \over \lambda^{\alpha}} \left[ 1 + O({1 \over \lambda}) \right]
~~\mbox{for $\lambda \vain + \infty$}~. \label{lapint} \feq
\begin{prop}
\textbf{Lemma.} Let
\beq L_{\delta}(\lambda) :=
\int_{0}^1 d t \, {(1 - t)^{\lambda} \over \sqrt{t} (3 + t)^{3 \lambda + \delta} }
\qquad \mbox{for \, $\delta \in \reali$, $\lambda \in (0,+\infty)$}~. \label{intelle} \feq
Then, for each $\delta \in \reali$,
\beq L_{\delta}(\lambda) = {1 + O(1/\lambda) \over 3^{3 \lambda + \delta}} \sqrt{\pi \over 2 \lambda}
\qquad \mbox{for $\lambda \vain + \infty$}~. \label{teselle} \feq
\end{prop}
\textbf{Proof.} We have $L_{\delta}(\lambda) = \int_{0}^1 d t\, \theta_{\delta}(t) e^{-\lambda \varphi(t)}$, where
\beq \theta_{\delta}(t) := {1 \over \sqrt{t} (3  + t)^{\delta}}~,
\qquad \varphi(t) := 3 \log(3 + t) - \log(1 - t)~. \feq
It is easily checked that
\parn
\vbox{
\beq \varphi'(t)= {2 (3-t) \over (1-t) (3 + t)} > 0 ~~~\mbox{for $t \in [0,1)$}~, \feq
$$ \varphi(0) = 3 \log 3~, \qquad \varphi(t)- \varphi(0) = 2 t + O(t^2) ~~\mbox{for $t \vain 0^{+}$}~, $$
$$ {\theta_{\delta}(t) \over \varphi'(t)} = {(\varphi(t) - \varphi(0))^{-1/2} \over \sqrt{2} \, 3^{\delta}}
[ 1 + O(\varphi(t)- \varphi(0)) ] \quad \mbox{for $t \vain 0^{+}$}~; $$}
so, application of \rref{lapint} yields the thesis \rref{teselle}. \fine
\parn
\textbf{Proof of Proposition \ref{co22}.} As usually, we consider any fixed space dimension
$d \in \naturali_0$. We must prove condition $\Co_{2 2 d}$, i.e.,
\beq \sup_{u \in [0,+\infty)} \Sigma_{2 2 d \ell}(u) = 1 + O(1/\ell) \qquad \mbox{for $\ell
\vain + \infty$}~. \label{co22d} \feq
Due to Eqs. \rref{bcell} \rref{eqsig}, for each $u \geqs 0$ we have
\beq \Sigma_{2 2 d \ell}(u) = (1 + 4 u)^\ell~
{}_{3} F_{2}(4 \ell - {d \over 2}, 2 \ell, 2 \ell;
2 \ell, 2 \ell + {1 \over 2}; -u) \label{dfu} \feq
$$ = (1 + 4 u)^\ell~
{}_{2} F_{1}(4 \ell - {d \over 2}, 2 \ell; 2 \ell + {1 \over 2}; -u)
\qquad \mbox{for $u \geqs 0$, $\ell > d/6$}~$$
(the last equality depends on Eq. \rref{ofco}). Now, using for
${}_{2} F_{1}$ the integral representation \rref{irep} we get
\beq \Sigma_{2 2 d \ell}(u) = {\Gamma(2 \ell + 1/2) \over \sqrt{\pi} \, \Gamma(2 \ell)}
\int_{0}^1 d s \, {s^{2 \ell - 1} \over \sqrt{1 - s}} \, W_{s d \ell}(u),\quad
W_{s d \ell}(u) := {(1 + 4 u)^\ell \over (1 + s u)^{4 \ell - d/2}}; \feq
of course, this implies
\beq \sup_{u \in [0,+\infty)}
\Sigma_{2 2 d \ell}(u) \leqs {\Gamma(2 \ell + 1/2) \over \sqrt{\pi} \, \Gamma(2 \ell)}
\int_{0}^1 d s \, {s^{2 \ell - 1} \over \sqrt{1 - s}} \left(\sup_{u \in [0,+\infty)} W_{s d \ell}(u)
\right)~. \label{ofcim} \feq
For all $\ell > d/6$ and $s \in (0,1)$, the function $W_{s d \ell} : [0,+\infty)
\vain (0, + \infty)$ attains its maximum at the point
\beq u_{s d \ell} := {1 - (1 - {d \over 8 \ell})s  \over 3 (1 - {d \over 6 \ell}) s}~, \feq
and so
\beq \sup_{u \in [0,+\infty)} W_{s d \ell}(u) = W_{s d \ell}(u_{s d \ell})
= {({3 \over 4})^{3 \ell - d/2} (1 - {d \over 6 \ell})^{3 \ell - d/2}
\over (1 - {d \over 8 \ell})^{4 \ell - d/2} s^\ell (1 - {s \over 4})^{3 \ell - d/2}}~. \feq
Inserting this result into Eq. \rref{ofcim} we get
\parn
\vbox{
\beq \sup_{u \in [0,+\infty)} \Sigma_{2 2 d \ell}(u) \leqs U_{d \ell}~, \label{summar1} \feq
$$ U_{d \ell} :=
{({3 \over 4})^{3 \ell - d/2} (1 - {d \over 6 \ell})^{3 \ell - d/2}
\over (1 - {d \over 8 \ell})^{4 \ell - d/2}}\,
{\Gamma(2 \ell + 1/2) \over \sqrt{\pi} \Gamma(2 \ell)}\,
\int_{0}^1 d s \, {s^{\ell - 1} \over \sqrt{1 - s} \, \, (1 - {s \over 4})^{3 \ell - d/2} }~; $$}
now, with a change of variable $s = 1-t$ in the integral and a comparison with Eq. \rref{intelle}, we find that
\beq U_{d \ell} = {3^{3 \ell - d/2} (1 - {d \over 6 \ell})^{3 \ell - d/2}
\over (1 - {d \over 8 \ell})^{4 \ell - d/2}}\,
{\Gamma(2 \ell + 1/2) \over \sqrt{\pi} \Gamma(2 \ell)}\, L_{3 - d/2}(\ell-1)~ \label{defuld} \feq
(the last factor indicates the Laplace integral $L_{\delta}(\lambda)$ of Eq.
\rref{intelle}, with $\lambda= \ell-1$ and $\delta = 3 - d/2$). Let us determine the
behavior of $U_{d \ell}$ for $\ell \vain + \infty$. To this purpose, we use the relations
\beq \big(1 - {d \over 6 \ell}\big)^{3 \ell - d/2} = e^{-d/2} [1 + O\big({1 / \ell}\big)],
\quad \big(1 - {d \over 8 \ell}\big)^{4 \ell - d/2} = e^{-d/2} \big[1 + O\big({1/\ell}\big) \big],
\label{relfour} \feq
$$ {\Gamma(2 \ell + 1/2) \over \Gamma(2 \ell)} = \sqrt{2 \ell} \big[1 + O\big({1/\ell}\big) \big],
\quad L_{3 - d/2}(\ell-1) = {1 + O(1/\ell) \over 3^{3 \ell - d/2}}
\sqrt{\pi \over 2 \ell}~; $$
the first two are obvious, the third one follows from Eq. \rref{rapgamma} and the fourth one
comes from the asymptotics \rref{teselle} of $L_{\delta}(\lambda)$. Inserting the relations
\rref{relfour} into \rref{defuld}, we get
\beq U_{d \ell} = 1 + O({1/\ell})~. \label{summar2} \feq
Let us summarize Eqs. \rref{summar1} \rref{summar2}:
\beq \sup_{u \in [0,+\infty)} \Sigma_{2 2 d \ell}(u) \leqs U_{d \ell} =  1 + O(1/\ell)
\qquad \mbox{for $\ell \vain + \infty$}~; \label{summar3} \feq
obviously enough, it is also
\beq \sup_{u \in [0,+\infty)} \Sigma_{2 2 d \ell}(u) \geqs \Sigma_{2 2 d \ell}(0)  = 1
\label{summar4} \feq
and Eqs. \rref{summar3} \rref{summar4} give the thesis \rref{co22d}. \fine
\begin{rema}
\label{remsig}
\textbf{Remark.}
Using Eq. \rref{dfu} with the known relation $(d / dw) \Big\vert_{w=0} \, {}_2 F_1(a, b, c, w) = a b/c$,
one easily finds that
\beq \left. {d \over d u} \right|_{u=0} \Sigma_{2 2 d \ell}(u) = {2 (d + 2) \ell \over 4 \ell + 1} > 0~. \feq
So, the function $\Sigma_{2 2 d \ell} : [0,+\infty) \vain (0,+\infty)$
is strictly increasing in a neighborhood of $u=0$; we also remark that
$\left. (d /d u) \right|_{u=0} \Sigma_{2 2 d \ell}(u) \vain {d / 2} + 1$ for
$\ell \vain + \infty$.
Even though  $u=0$ is not a maximum point, the $\ell \vain + \infty$
asymptotics $\sup_{u \geqs 0} \Sigma_{2 2 d \ell}(u) = 1 + O(1/\ell)=
\Sigma_{2 2 d \ell}(0) + O(1/\ell)$ suggests that, for large $\ell$,
the sup of $\Sigma_{2 2 d \ell}$ could be obtained at a point $O(1/\ell)$. We have
found numerical evidence for this: $\Sigma_{2 2 d \ell}$ seems to have a unique maximum point
$u_{2 2 d \ell}$, such that
$u_{2 2 d \ell} = O(1/\ell)$ for $\ell \vain + \infty$. \fine
\end{rema}
\begin{prop}
\textbf{Lemma.} Let
$f_{\sigma d}(x) := e^{-\si |x|^2/2}$ for $x \in \reali^d$ and $\si >0$, as in \rref{fsigd};
furthermore, fix
\beq a \in (0,+\infty), \qquad \zeta \in (0,1/a)~. \label{azeta} \feq
Then, with $\|~\|_{a \ell}$ indicating the $H^{a \ell}$ norm,
\beq \| f_{\zeta/\ell, d} \|_{a \ell} =
\left[{\pi \ell \over \zeta (1 - a \zeta)}\right]^{d/4} \left[1 + O({1 \over \ell})\right]  \qquad \mbox{for $\ell
\vain + \infty$}~. \label{givebyas}\feq
\end{prop}
\textbf{Proof. } Eq. \rref{giveby0} gives
\beq \| f_{\zeta/\ell, d} \|_{a \ell}^2 = {2 \, \pi^{d/2} \ell^d \over \Gamma(d/2) \zeta^{d}}~
\int_{0}^{+\infty} d \rho \,\rho^{d-1} (1 + \rho^2)^{a \ell} e^{-{\ell \rho^2 \over \zeta}}~; \feq
with a change of variable $\rho = \sqrt{\zeta t}$, we get
\beq \| f_{\zeta/\ell, d} \|_{a \ell}^2 = {\pi^{d/2} \ell^d \over \Gamma(d/2) \zeta^{d/2}} L_{a \zeta d}(\ell)~,
\qquad L_{a \zeta d}(\ell) := \int_{0}^{+\infty} d t \, t^{d/2-1} (1 + \zeta t)^{a \ell} e^{-\ell t}~.
\label{gibeby1} \feq
We note that
\beq L_{a \zeta d}(\ell) = \int_{0}^{+\infty} d t \, \vartheta_{d}(t) e^{\dd{-\ell \varphi_{a \zeta}(t)}}~, \feq
$$ \vartheta_{d}(t) := t^{d/2-1}~, \qquad \varphi_{a \zeta}(t) := t - a \log(1 + \zeta t)~; $$
this indicates that $L_{a \zeta}(\ell)$ is a Laplace integral in the parameter $\ell$, in the sense
reviewed at the beginning of the section. One easily checks that
\beq {\varphi'}_{a \zeta}(t) =
{1 - a \zeta + \zeta t \over 1 + \zeta t} > 0 \quad \mbox{for $t \in [0,+\infty)$}~, \feq
$$ \varphi_{a \zeta}(0) = 0~, \qquad
\varphi_{a \zeta}(t) = (1 - a \zeta) t + O(t^2) \quad \mbox{for $t \vain 0^{+}$}~, $$
$$ {\vartheta_d(t) \over {\varphi'}_{a \zeta}(t)} = {\varphi_{a \zeta}(t)^{d/2-1} \over (1 - a \zeta)^{d/2}}\,
 [1 + O(\varphi_{a \zeta}(t))]
\quad \mbox{for $t \vain 0^{+}$}~; $$
from here and \rref{lapint}, we get
\beq L_{a \zeta d}(\ell) = {\Gamma(d/2) \over (1 - a \zeta)^{d/2} \, \ell^{d/2}}
\left[1 + O({1 \over \ell})\right] \qquad \mbox{for $\ell \vain + \infty$}~. \label{insel} \feq
Inserting \rref{insel} into \rref{gibeby1}, and taking the square root, we get the thesis
\rref{givebyas}. \fine
\textbf{Proof of Proposition \ref{kfas}.} Let $1 \leqs b \leqs c$ and
$\xi \in (0,1/b)$, $\eta \in (0,1/c)$ with $\xi + \eta < 1$; we must derive the
$\ell \vain + \infty$ asymptotics
\rref{relk}, i.e.,
\beq { \| f_{\xi/\ell + \eta/\ell, d} \|_\ell \over \| f_{\xi/\ell,  d} \|_{b \ell} \| f_{\eta/\ell, d}\|_{c \ell}}
= {1 + O(1/\ell) \over [D_{b c}(\xi,\eta) \pi \ell]^{d/4}}~, \qquad
D_{b c}(\xi,\eta) := {(1 - \xi - \eta) (\xi + \eta) \over \xi \eta (1 - b \xi) (1 - c \eta)}~. \label{relkk} \feq
The thesis follows using Eq. \rref{givebyas} with $(a, \zeta) = (1, \xi+ \eta)$, or $(b, \xi)$, or $(c, \eta)$
(in each of the three cases, the assumptions on $\xi, \eta$ ensure conditions \rref{azeta} to be fulfilled).
\fine
\vskip 0.5cm \noindent
\textbf{Acknowledgments.}
This work has been partially supported by the GNFM
of Istituto Nazionale di Alta Ma\-te\-ma\-ti\-ca and by MIUR,
Research Project Cofin/2006
``Geometric methods in the theory of nonlinear waves and their applications''.
\vskip 0.5cm \noindent
\appendix
\section{Appendix. Derivation of Eq.\rref{result}.}
\label{appef}
Let us consider the integral
\beq I_{\mu \nu \de}(h) := \int_{0}^{+\infty} \!\!\! dr\, r^{\mu + \nu + \delta + 1} J_{\delta}(h r) K_{\mu}(r) K_{\nu}(r)~;
\label{imunu} \feq
with this notation, Eq. \rref{result} reads
\parn
\vbox{
\beq I_{\mu \nu \de}(h) = 2^{\mu + \nu + \delta - 1}~
{\Gamma(\mu + \delta + 1) \Gamma(\nu + \delta + 1) \Gamma(\mu + \nu + \delta + 1) \over
\Gamma(\mu + \nu + 2 \delta + 2 )} \, h^{\delta} \label{resultt} \feq
$$ \times \,
{}_{3} F_{2}(\mu + \delta + 1, \nu + \delta + 1, \mu + \nu + \delta + 1;
{\mu + \nu \over 2} + \delta + 1, {\mu + \nu \over 2} + \delta + {3 \over 2}; -{h^2 \over 4}) $$
$$ \mbox{for $h,\mu,\nu,\delta \in \reali$,~ $h > 0$,
$\delta, \mu + \delta, \nu + \delta, \mu + \nu + \delta > -1$}~. $$}
In the sequel we prove this identity, after checking preliminarily that the integral in the
right hand side converges under the above conditions for $h, \mu, \nu, \delta$. \parn
Convergence of the integral follows immediately from the relations
$J_{\xi}(w) = O(w^{\xi})$, $K_{\eta}(w) = O(w^{-|\eta|})$ for $\xi > -1$,
$\eta \in \reali$, $w \vain 0^{+}$ and $J_{\xi}(w) = O(1/\sqrt{w})$,
$K_{\eta}(w) = e^{-w} O(1/\sqrt{w})$ for $\xi, \eta \in \reali$, $w \vain +\infty$
(see \cite{Wat}, Chapters III and VII);
these ensure integrability of the function of $r$ in $I_{\mu \nu \delta}(h)$, both
near zero and near $+\infty$.
To derive the equality \rref{resultt},
we start from the familiar series expansion (see again \cite{Wat}, Chapter III)
\beq J_{\delta}(w) = \sum_{k=0}^{+\infty} {(-1)^k \over k! \Gamma(\de + 1 + k)} ({w \over 2})^{\de + 2 k}~,
\feq
to be applied with $w = h r$; inserting this into Eq. \rref{imunu}, we get
\beq I_{\mu \nu \de}(h)
= ({h \over 2})^{\de}
\sum_{k=0}^{+\infty} {1 \over k! \Gamma(\de + 1 + k)} ({-h^2 \over 4})^{k}
\int_{0}^{+\infty} \!\!\! dr\, r^{2 \delta + \mu + \nu + 1 + 2 k} K_{\mu}(r) K_{\nu}(r)~.
\label{kapkap} \feq
On the other hand,
\beq \int_{0}^{+\infty} \!\!\! dr\, r^{\alpha-1} K_{\mu}(r) K_{\nu}(r) \label{useq} \feq
$$ = {2^{\al - 3} \over \Gamma(\alpha)} \Gamma({\al - \mu - \nu\over 2})
\Gamma({\al + \mu - \nu\over 2}) \Gamma({\al - \mu + \nu\over 2})
\Gamma({\al + \mu + \nu\over 2}) $$
if the arguments of all the above Gamma functions are positive
(this is a special case of an identity in \cite{Gra}: see Eq. (6.576.4), page 693).
We can use Eq. \rref{useq} to compute the integrals in \rref{kapkap}, the conclusion being
\par \noindent
\vbox{
\beq I_{\mu \nu \de}(h) \feq
$$ = 2^{\mu + \nu + \de - 1} h^{\de}
\sum_{k=0}^{+\infty} {1 \over k!} {\Gamma(\mu + \de + 1 + k)
\Gamma(\nu + \de + 1 + k) \Gamma(\mu + \nu + \de + 1 + k) \over \Gamma(\mu + \nu + 2 \de + 2 + 2 k)}
(-h^2)^{k}~. $$ }
Now, we introduce the relations
\beq \Gamma(\alpha + k) = (\alpha)_k \Gamma(\alpha), \quad
\Gamma(2 \alpha + 2 k) = 4^{k} (\alpha)_k (\alpha + {1 \over 2})_k \Gamma(2 \alpha)
\quad \mbox{for $k \in \naturali$} \feq
(the first appearing in Eq.\rref{elem}, the second following from the first and from the elementary identity
$(2 \alpha)_{2 k} = 4^{k} (\alpha)_k (\alpha + 1/2)_k$). In this way we get
\beq I_{\mu \nu \de}(h) = 2^{\mu + \nu + \de - 1}
{\Gamma(\mu + \de + 1)
\Gamma(\nu + \de + 1) \Gamma(\mu + \nu + \de + 1) \over \Gamma(\mu + \nu + 2 \de + 2)}
h^{\de}~\feq
$$ \times\, \sum_{k=0}^{+\infty} {1 \over k!} {(\mu + \de + 1)_k
(\nu + \de + 1)_k (\mu + \nu + + \de + 1)_k \over ({\mu + \nu \over 2} + \de + 1)_k
({\mu + \nu \over 2} + \de + {3 \over 2})_{k}} (-{h^2 \over 4})^{k}~. $$
According to Eq. \rref{partpq}, the above series equals
$${}_{3} F_{2}(\mu + \delta + 1, \nu + \delta + 1, \mu + \nu + \delta + 1;
{\mu + \nu \over 2} + \delta + 1, {\mu + \nu \over 2} + \delta + {3 \over 2}; -{h^2 \over 4})~,$$
so Eq. \rref{resultt} is proved. (Final remark: in fact, the previous considerations
give the thesis \rref{resultt} for $h^2/4 < 1$, i.e. $h \in (0,2)$, since the series expansion \rref{partpq}
for ${}_{3} F_{2}$ has a convergence radius $1$. However, after proving the thesis
for $h \in (0,2)$ one can extend it to all $h \in (0,+\infty)$ by a standard
application of the analytic continuation principle.)
\section{Appendix. Calculation of the upper and lower bounds $\boma{K^{\pm}_{\ell m n d}}$
in the table of page \pageref{patable}: some examples.}
\label{apptab}
\textbf{Computation of $\boma{K^{+}_{0 1 2 1}}$.}
(i) We first determine the $\SS$-function upper bound. Eqs.(\ref{ff}-\ref{spe2}) give
\beq \SS_{0 1 2 1}(u) = {3 + u \over 16(1 + u)^2}~\qquad \mbox{for $u \in [0,+\infty)$}~; \feq
the above function is easily studied by analytical means, the conclusion being
\beq \sup_{u \in [0,+\infty)}  \SS_{0 1 2 1}(u) = \SS_{0 1 2 1}(0) = {3 \over 16}~. \feq
So,
\beq \sqrt{\sup_{u \in [0,+\infty)}  \SS_{0 1 2 1}(u)} = \sqrt{{3 \over 16}} \leqs 0.434
:= K^{\SS}_{0 1 2 1}~.  \label{434} \feq
(ii) Let us build the H\"older upper bound \rref{kimbsuu}; in this case,
Eqs. \rref{pert} \rref{kimbsu0}
\rref{kimbsu} give $R_{1 1} = R_{2 1} = R_{1 2 1} = [2, +\infty]$, so we must evaluate
$\inf_{p \in [2, +\infty]} S_{p 1 1} \, S_{p^* 2 1}$,
the factors $S_{p 1 1}$, $S_{p^* 2 1}$ being given by Eqs. (\ref{srid}-\ref{siid}).
As found numerically, the inf is attained for $p$ close to $3.21$, and
\beq \inf_{p \in [2, +\infty]} S_{p 1 1} \, S_{p^* 2 1} \leqs 0.383 := K^{\Ii}_{0 1 2 1}~. \feq
(iii) The H\"older bound  $K^{\Ii}_{0 1 2 1}$ is better than the
$\SS$-function bound $K^{\SS}_{0 1 2 1}$, so we take
\beq K^{+}_{0 1 2 1} := K^{\Ii}_{0 1 2 1} = 0.383~; \label{ki} \feq
this is the value reported in the table.
\vskip 0.2cm \noindent
\textbf{Computation of $\boma{K^{-}_{0 1 2 1}}$.}
(i) We first consider the Bessel lower bound
\rref{kbes} with $s = 1$, $t=2$. In this case, Eqs.
\rref{fnutd} \rref{pooint}  give
\beq \| \f_{\mu 1 1} \|^2_{1} = {\pi \over 2} \, {1 + \mu^2 \over \mu}~, \qquad
\| \f_{\nu 2 1} \|^2_{2} = {\pi \over 16} \, {5  + 2 \nu^2 + \nu^4 \over \nu}~, \feq
\beq \| \f_{\mu 1 1} \f_{\nu 2 1} \|^2_0 = {\pi^2 \over 32} \,
{2 \mu^2 + 6 \mu \nu + 5 \nu^2 \over (\mu + \nu)^3}; \feq
from here one computes, according to Eq. \rref{theabove}, the function
\beq
\Kap^{B_{1 2}}_{0 1 2 1}(\mu, \nu) :=
{ \| \f_{\mu 1 1} \f_{\nu 2 1} \|_0 \over \| \f_{\mu 1 1} \|_1 \| \f_{\nu 2 1}\|_2}~
\qquad (\mu, \nu \in (0,+\infty))~. \feq
It is found numerically that the above function attains its sup for $(\mu, \nu)$ close to
$(0.499, 0.784)$, and that
\beq \sup_{\mu, \nu > 0} \Kap^{B_{1 2}}_{0 1 2 1}(\mu, \nu)
\geqs 0.951 \, K^{+}_{0 1 2 1} := K^{B_{ 1 2}}_{0 1 2 1}~. \label{pb} \feq
(ii) We pass to the Fourier lower bound \rref{thea111}. In this case,
from Eq. \rref{gv} one gets
\beq \| f_{h \kappa 1} \|^2_ 0 = \sqrt{{\pi \over \kappa}}~, \qquad
\| f_{p \si 1} \|^2_ 1 = {\sqrt{\pi} \over 2} \, {2 + 2 p^2 + \si \over \sqrt{\si}}~, \feq
$$ \| f_{q \tau 1} \|^2_ 2 = {\sqrt{\pi} \over 4} \, {4 + 8 q^2 + 4 q^4 + 4 \tau + 12 q^2 \tau +
3 \tau^2 \over \sqrt{\tau}} $$
for $h,p,q \in [0,+\infty)$ and $\kappa, \si,\tau \in (0,+\infty)$;
from here, one computes the function
\beq \Kap^{F}_{0 1 2 1}(p,q,\si,\tau) := {\| f_{p+q, \si+\tau, 1} \|_ 0 \over
\| f_{p \si 1} \|_1 \, \| f_{q \tau 1} \|_2} \feq
($p, q \in [0,+\infty)$, $\si, \tau \in (0,+\infty)$). A numerical investigation
seems to indicate that the sup of this function is attained
for $(p, q, \si, \tau)$ close to $(0,0, 0.472, 0.291)$; in any case, using the value
at this point as a lower approximant for the sup we get
\beq
\sup_{p, q \geqs 0, \, \si, \tau > 0} \Kap^{F}_{0 1 2 1}(p, q, \si,\tau) \geqs 0.987 \, K^{+}_{0 1 2 1}
:= K^{F}_{0 1 2 1}~. \feq
(iii) The Fourier lower bound $K^{F}_{0 1 2 1}$ is better than the Bessel lower
bound $K^{B_{1 2}}_{0 1 2 1}$; in conclusion we take
\beq K^{-}_{0 1 2 1} :=   K^{F}_{0 1 2 1} = 0.987 \, K^{+}_{0 1 2 1}~, \feq
as indicated in the table. The symbol $(F)$ appearing in the table
recalls that the lower bound $K^{-}_{0 1 2 1}$ is of the Fourier type.
\vskip 0.2cm \noindent
\textbf{Computation of $\boma{K^{+}_{4 5 6 1}}$.} We use for this the $\SS$-function upper
bound. Eqs. (\ref{ff}-\ref{spe2}) give
\beq \SS_{4 5 6 1}(u) =
(1 + 4 u)^4 \, {46189 + 20995 \, u + 9690 \, u^2 + 3230 \, u^3 + 665 \, u^4 + 63 \, u^5 \over
 524288 \, (1 + u)^{10}}~ \feq
for $u \in [0,+\infty)$. It is found numerically that the above function attains its sup
close to $u = 0.315$, and that
\beq \sqrt{\sup_{u \in [0,+\infty)}  \SS_{4 5 6 1}(u)} \leqs 0.417 :=
K^{+}_{4 5 6 1}~; \feq
this upper bound is reported in the table.
\vskip 0.2cm \noindent
\textbf{Computation of $\boma{K^{-}_{4 5 6 1}}$.}
(i) We first consider the Bessel lower bound
\rref{kbes} with $s = 5$, $t=6$. Eq.
\rref{fnutd}   gives
\beq \| \f_{\mu 5 1} \|^2_{5} =
{5 \pi \over 65536 \mu} \, (2431 + 715 \mu^2 + 286 \mu^4 + 110 \mu^6 + 35 \mu^8 + 7 \mu^{10})~, \feq
$$ \| \f_{\nu 6 1} \|^2_{6} =
{3 \pi \over 524288 \nu} (29393 + 8398 \nu^2 + 3315 \nu^4 + 1300 \nu^6 + 455 \nu^8 +
126 \nu^{10} + 21 \nu^{12})~$$
for $\mu, \nu \in (0,+\infty)$.
Eq. \rref{pooint} gives
\beq \| \f_{\mu 5 1} \f_{\nu 6 1} \|^2_4 =
{\pi^2 \over 34359738368 \, (\mu + \nu)^{19}} \, P(\mu,\nu) \feq
where $P(\mu, \nu)$ is a polynomial of the form:
\beq P(\mu, \nu) = \sum_{i, j \in \naturali, 18 \leqs i+j \leqs 26} P_{i j} \mu^i \nu^j~,
\qquad P_{i j} \in \naturali~\mbox{for all $i, j$}~. \feq
The full expression of this polynomial is easily computed with MATHEMATICA, but it is too
long to be reported here; as examples we give only three coefficients, namely,
\beq P_{18, 0} = 192972780, \qquad P_{1, 25} = 4236050, \qquad P_{0, 26} = 222950. \feq
The expressions of the above norms determine the function
\beq
\Kap^{B_{5 6}}_{4 5 6 1}(\mu, \nu) :=
{ \| \f_{\mu 5 1} \f_{\nu 6 1} \|_4 \over \| \f_{\mu 5 1} \|_5 \| \f_{\nu 6 1}\|_6}~
\qquad (\mu, \nu \in (0,+\infty))~. \feq
It is found numerically that the above function attains its sup for $(\mu, \nu)$ close to
$(1.19,1.14 )$, and that
\beq \sup_{\mu, \nu > 0} \Kap^{B_{5 6}}_{4 5 6 1}(\mu, \nu)
\geqs 0.823 \, K^{+}_{4 5 6 1} := \Kap^{B_{5 6}}_{4 5 6 1}~. \label{pbb} \feq
(ii) We pass to the Fourier lower bound \rref{thea111}.
From Eq. \rref{gv} one gets
\beq \| f_{h \kappa 1} \|^2_ 4 =
{1 \over 16} \sqrt{\pi \over \kappa} \,
(16 + 64 h^2 + 96 h^4 + 64 h^6 + 16 h^8 + 32 \kappa + 288 h^2 \kappa \feq
$$ + 480 h^4 \kappa + 224 h^6 \kappa + 72 \kappa^2 + 720 h^2 \kappa^2 + 840 h^4 \kappa^2 +
   120 \kappa^3 + 840 h^2 \kappa^3 + 105 \kappa^4)~, $$
$$ \| f_{p \si 1} \|^2_ 5 = {1 \over 32} \sqrt{\pi \over \si} \,
(32 + 160 p^2 + 320 p^4 + 320 p^6 + 160 p^8 + 32 p^{10} + 80 \si $$
$$ +   960 p^2 \si + 2400 p^4 \si + 2240 p^6 \si + 720 p^8 \si + 240 \si^2 +
   3600 p^2 \si^2 + 8400 p^4 \si^2 $$
$$ + 5040 p^6 \si^2 + 600 \si^3 + 8400 p^2 \si^3 +
12600 p^4 \si^3 + 1050 \si^4 + 9450 p^2 \si^4 + 945 \si^5)~, $$
$$ \| f_{q \tau 1} \|^2_ 6 =  {1 \over 64} \sqrt{\pi \over \tau} (
64 + 384 q^2 + 960 q^4 + 1280 q^6 + 960 q^8 + 384 q^{10} + 64 q^{12} + 192 \tau $$
$$ + 2880 q^2 \tau + 9600 q^4 \tau + 13440 q^6 \tau + 8640 q^8 \tau +
2112 q^{10} \tau + 720 \tau^2 + 14400 q^2 \tau^2 $$
$$  + 50400 q^4 \tau^2 +
60480 q^6 \tau^2 + 23760 q^8 \tau^2 + 2400 \tau^3 + 50400 q^2 \tau^3 + 151200 q^4 \tau^3 $$
$$  + 110880 q^6 \tau^3 + 6300 \tau^4 + 113400 q^2 \tau^4 +
207900 q^4 \tau^4 + 11340 \tau^5 + 124740 q^2 \tau^5 + 10395 \tau^6) $$
for $h, p, q \in [0,+\infty)$ and $\kappa, \si,\tau \in (0,+\infty)$;
from here, one computes the function
\beq \Kap^{F}_{4 5 6 1}(p,q,\si,\tau) := {\| f_{p+q, \si+\tau, 1} \|_ 4 \over
\| f_{p \si 1} \|_5 \, \| f_{q \tau 1} \|_6} \feq
($p, q \in [0,+\infty)$, $\si, \tau \in (0,+\infty)$). A numerical investigation
seems to indicate that the sup of this function is attained
for $(p, q, \si, \tau)$ close to $( 0.288,$ $0.215,$ $0.147,$ $0.109)$; in any case,
using the value at this point as a lower approximant for the sup we get
\beq
\sup_{p, q \geqs 0, \si, \tau > 0} \Kap^{F}_{4 5 6 1}(p,q,\si, \tau)
\geqs 0.878 \, K^{+}_{4 5 6 1}
:= K^{F}_{4 5 6 1}~. \feq
(iii) The Fourier lower bound $K^{F}_{4 5 6 1}$ is better than the Bessel lower
bound $K^{B_{5 6}}_{4 5 6 1}$; in conclusion we take
\beq K^{-}_{4 5 6 1} :=  K^{F}_{4 5 6 1} = 0.878 \, K^{+}_{4 5 6 1}~, \feq
as indicated in the table. The symbol $(F)$ appearing in the table
recalls the type of the lower bound $K^{-}_{4 5 6 1}$.
\vskip 0.2cm \noindent
\textbf{Computation of $\boma{K^{+}_{1 1 2 3}}$.}
We use for this the $\SS$-function upper
bound. Eqs.(\ref{ff}-\ref{spe2}) give
\beq \SS_{1 1 2 3}(u) = {(1 + 4 u) \over 32 \pi (1 + u)}~ \feq
for $u \in [0,+\infty)$. The above function attains its sup in the limit $u
\vain +\infty$, and
\beq \sqrt{\sup_{u \in [0,+\infty)}  \SS_{1 1 2 3}(u)} = \sqrt{\SS_{1 1  2 3}(+\infty)} =
{1 \over 2 \sqrt{2 \pi}} := K^{+}_{1 1 2 3}~. \label{eqqa} \feq
This is the value reported in the table; from a numerical viewpoint, $K^{+}_{1 1 2 3} =  0.1994...$~.
\vskip 0.2cm \noindent
\textbf{Computation of $\boma{K^{-}_{1 1 2 3}}$.}
We are discussing a case with $\ell=m$, so we have the $S$-constant lower bound \rref{kell}; more precisely, this
bound is (recalling Eq. \rref{siid})
\beq S_{\infty 2 3} = {1 \over 2 \sqrt{2 \pi}} := K^{-}_{1 1 2 3}~. \label{eqqb} \feq
This lower bound equals $K^{+}_{1 1 2 3}$; we can avoid calculating the
Bessel and Fourier lower bounds, since they cannot be better. In the table
we have indicated that $K^{-}_{1 1 2 3}/K^{+}_{1 1 2 3}=1$, and we have used the symbol (S)
to recall the type of the lower bound.
\parn
Of course, in this case we have the sharp constant:
\beq K_{1 1 2 3} = K^{\pm}_{1 1 2 3}~\label{eqqc}~. \feq
\vskip 0.2cm \noindent
\textbf{Computation of $\boma{K^{+}_{2 2 3 3}}$.} Again, we use the $\SS$-function bound.
Eqs.(\ref{ff}-\ref{spe2}) give
\beq \SS_{2 2 3 3}(u) = {(1 + 4 u)^2 (5 + u) \over 512 \pi (1 + u)^3}~ \feq
for $u \in [0,+\infty)$. It is found that
\beq \sqrt{\sup_{u \in [0,+\infty)}  \SS_{2 2 3 3}(u)} = \sqrt{\SS_{2 2  3 3}({{13 \over 5}})} =
{19 \over 288} \sqrt{{19 \over 2 \pi}} \leqs 0.115 := K^{+}_{2 2 3 3}~; \feq
this upper bound is reported in the table.
\vskip 0.2cm \noindent
\textbf{Computation of $\boma{K^{-}_{2 2 3 3}}$.}
(i) Let us compute the Bessel lower bound
\rref{kbes}, with $s = 2$, $t=3$. Eqs.
\rref{fnutd} \rref{pooint}  give
\beq \| \f_{\mu 2 3} \|^2_{2} = {\pi^2 \over 8 \mu^3} (1 + 2 \mu^2 + 5 \mu^4)~, \qquad
\| \f_{\nu 3 3} \|^2_{3} = {\pi^2 \over 128 \nu^3} \, (7 + 9 \nu^2 + 9 \nu^4 + 7 \nu^6)~, \feq
\beq \| \f_{\mu 2 3} \f_{\nu 3 3} \|^2_2 = {\pi^3 \over 1024 (\mu + \nu)^5} \,
(\mu^2 + 2 \mu^4 + 5 \mu^6 + 5 \mu \nu + 10 \mu^3 \nu + 25 \mu^5 \nu \feq
$$ + 7 \nu^2 +
   20 \mu^2 \nu^2 + 53 \mu^4 \nu^2 + 18 \mu \nu^3 + 62 \mu^3 \nu^3 + 6 \nu^4 +
   43 \mu^2 \nu^4 + 17 \mu \nu^5 + 3 \nu^6)~; $$
from here one computes, according to Eq. \rref{theabove}, the function
\beq
\Kap^{B_{2 3}}_{2 2 3 1}(\mu, \nu) :=
{ \| \f_{\mu 2 3} \f_{\nu 3 3} \|_2 \over \| \f_{\mu 2 3} \|_2 \| \f_{\nu 3 3}\|_3}~
\qquad (\mu, \nu \in (0,+\infty))~. \feq
It is found numerically that the above function attains its sup for $(\mu, \nu)$ close to
$(1.31, 1.04)$, and that
\beq \sup_{\mu, \nu > 0} \Kap^{B_{ 2 3}}_{2 2 3 1}(\mu, \nu)
\geqs 0.916 \, K^{+}_{2 2 3 1} := \Kap^{B_{ 2 3}}_{2 2 3 1}~. \label{pbc} \feq
(ii) Let us pass to the Fourier lower bound \rref{thea111}.
From Eq. \rref{gv} one gets
\parn
\vbox{
\beq
\| f_{p \si 3} \|^2_ 2 = {1 \over 4} \left( {\pi \over \si}\right)^{3/2} \!\!
(4 + 8 p^2 + 4 p^4 + 12 \si + 20 p^2 \si + 15 \si^2)~ , \feq
$$ \| f_{q \tau 3} \|^2_ 3 = {1 \over 8} \left( {\pi \over \tau}\right)^{3/2} \!\!
(8 + 24 q^2 + 24 q^4 + 8 q^6 + 36 \tau + 120 q^2 \tau + 84 q^4 \tau +
   90 \tau^2 + 210 q^2 \tau^2 + 105 \tau^3)~, $$}
for $p,q \in [0,+\infty)$ and $\si,\tau \in (0,+\infty)$;
from here, one computes the function
\beq \Kap^{F}_{2 2 3 3}(p,q,\si,\tau) := {\| f_{p+q, \si+\tau, 3} \|_ 2 \over
\| f_{p \si 3} \|_2 \, \| f_{q \tau 3} \|_3} \feq
($p, q \in [0,+\infty)$, $\si, \tau \in (0,+\infty)$). A numerical investigation
seems to indicate that the sup of this function is attained for
$(p, q, \si, \tau)$ close to $(0.667,$ $0.114,$ $2.53,$ $0.430)$; in any case,
using the value at this point as a lower approximant for the sup we get
\beq \sup_{p, q \geqs 0, \si, \tau > 0} \Kap^{F}_{2 2 3 3}(p,q,\si,\tau)
 \geqs 0.809 \, K^{+}_{2 2 3 3} := K^{F}_{2 2 3 3}~. \feq
(iii) Since we are discussing a case with $\ell=m$, we have also the $S$-constant lower bound \rref{kell}; this
bound is (recalling Eq. \rref{siid})
\beq S_{\infty 3 3} = {1 \over 4 \sqrt{2 \pi}} = 0.8672... \, K^{+}_{2 2 3 3}~. \feq
(iv) The Bessel lower bound $K^{B_{ 2 3}}_{2 2 3 3}$ is better than the $S$-constant and Fourier lower
bounds $S_{\infty 3 3}$, $K^{F}_{2 2 3 3}$; in conclusion we take
\beq K^{-}_{2 2 3 3} :=   K^{B_{ 2 3}}_{2 2 3 3} = 0.916 \, K^{+}_{2 2 3 3}~, \feq
as indicated in the table. The symbol $(B_{ 2 3})$ appearing in the table
recalls the type of the lower bound.
\vfill \eject 
\end{document}